\documentclass[aap,MSNbibl,seceqn,nameyear,dvips]{arximspdf}

\doi{10.1214/10-AAP755}
\volume{21}
\issue{6}
\pubyear{2011}
\firstpage{2191}
\lastpage{2225}

\makeatletter
\newcommand{\R}{\mathbb{R}}
\newcommand{\xib}{\mathbf{\xi}}
\newcommand{\ab}{\mathbf{a}}
\newcommand{\theb}{\bolds{\theta}}

\newproclaim{remark}{Remark}
\newtheorem{prop}{Proposition}
\newtheorem{lemma}{Lemma}
\newproclaim{ass}{Assumption}

\makeatother

\begin{document}
\begin{frontmatter}

\title{Optimal arbitrage under model uncertainty}
\runtitle{Optimal arbitrage under model uncertainty}

\begin{aug}
\author[A]{\fnms{Daniel} \snm{Fernholz}\ead[label=e1]{df@danielfernholz.com}}
and
\author[B]{\fnms{Ioannis} \snm{Karatzas}\corref{}\thanksref{aut2,aut1}\ead[label=e2]{ik@enhanced.com}\ead[label=e3]{ik@math.columbia.edu}}
\thankstext{aut2}{This author is on partial leave
from the Department of Mathematics,
Columbia University, New York, New York 10027, USA (\printead{e3}).}
\thankstext{aut1}{Supported by NSF Grant DMS-09-05754.}
\runauthor{D. Fernholz and I. Karatzas}
\affiliation{Daniel Fernholz LLC and INTECH Investment Management}
\dedicated{Dedicated to Professor Mark H. A. Davis on the occasion of
his 65th birthday}
\address[A]{Daniel Fernholz LLC\\
1108 Lavaca, Suite 110-440\\
Austin, Texas 78701\\
USA\\
\printead{e1}} 
\address[B]{INTECH Investment Management\\
One Palmer Square, Suite 441 \\
Princeton, New Jersey 08542\\
USA\\
\printead{e2}}
\end{aug}

\received{\smonth{5} \syear{2010}}
\revised{\smonth{12} \syear{2010}}

%
\begin{abstract}
In an equity market model with ``Knightian'' uncertainty regarding
the relative risk and covariance structure of its assets, we
characterize in several ways the highest return relative to the market
that can be achieved using nonanticipative investment rules over a
given time horizon, and under any admissible configuration of model
parameters that might materialize. One characterization is in terms of
the smallest positive supersolution to a fully nonlinear parabolic
partial differential equation of the Hamilton--Jacobi--Bellman type.
Under appropriate conditions, this smallest supersolution is the
value function of an associated stochastic control problem, namely, the
maximal probability with which an auxiliary multidimensional diffusion
process, controlled in a manner which affects both its drift and
covariance structures, stays in the interior of the positive orthant
through the end of the time-horizon. This value function is also
characterized in terms of a stochastic game, and can be used to
generate an investment rule that realizes such best possible
outperformance of the market.
\end{abstract}

%

\end{frontmatter}

\section{Introduction}
\label{intro}

Consider an equity market with asset capitalizations $ \mathfrak{X} (t)
= ( X_1 (t), \ldots, X_n (t) )' \in(0, \infty)^n $ at time $ t \in[0,
\infty) $, and with covariance and relative risk rates $ \alpha(t,
\mathfrak{X}) = \{ \alpha_{ij} (t, \mathfrak{X}) \}_{1 \le i,j \le
n} $
and $ \vartheta(t, \mathfrak{X}) = ( \vartheta_1 (t, \mathfrak{X}),
\ldots,\allowbreak\vartheta_n (t, \mathfrak{X}))' $, respectively. At any
given time $ t $, these rates are \textit{nonanticipative functionals} of
past-and-present capitalizations $ \mathfrak{X} (s)$, $0 \le s \le t $;
they are not specified with precision but are, rather, subject to
``Knightian uncertainty.'' To wit, for a given collection
%
\begin{equation}
\label{1.b}
\mathbb{K} = \{ \mathcal{K}
(\mathbf{y}) \}_{\mathbf{y} \in\mathfrak{S}_n } ,\qquad  \mathfrak{S}_n :=
[0,\infty)^n \setminus\{\mathbf{0}\}
\end{equation}
of nonempty compact and convex subsets on $ \mathbb{R}^n \times
\mathbb{S}^n $, where $ \mathbb{S}^n $ is the space of real, symmetric,
positive definite $ (n \times n) $ matrices, and $ \mathbf{ 0} $ is the
origin in~$ \mathbb{R}^n $, they are subject to the constraint
%
\begin{equation}
\label{E.3}
(\vartheta(t, \mathfrak{X}), \alpha(t, \mathfrak{X}) ) \in\mathcal
{K} (\mathfrak{X} (t)) \qquad  \mbox{for all } t \in[0, \infty) .
\end{equation}
In other words, the pair $ (\vartheta, \alpha) $ must take values at
time $ t $ inside the compact, convex set $ \mathcal{K} (\mathfrak{X}
(t)) $ which is determined by the current location of the asset
capitalization process; but within this range, the actual value $
(\vartheta(t, \mathfrak{X}), \alpha(t, \mathfrak{X}) ) $ is allowed
to depend on past capitalizations as well. [To put it a little
differently: the constraint (\ref{E.3}) is not necessarily
``Markovian,'' as long as the sets in (\ref{1.b}) are not singletons.]

Under these circumstances, what is the highest return on investment
relative to the market that can be achieved using nonanticipative
investment rules, and with probability one under all possible market
model configurations that satisfy the constraints of (\ref{E.3})? What
are the weights in the various assets of an investment rule that
accomplishes this?

Answers: Subject to appropriate conditions, $1 / U (T, \mathfrak{X} (0)
) $ and
%
\begin{eqnarray}\label{FG}
X_i (t) D_i \log U \bigl( T-t, \mathfrak{X}(t) \bigr) + { X_i
(t) \over X_1(t)+\cdots+ X_n (t) } ,\nonumber\\[-8pt]\\[-8pt]
\eqntext{i = 1, \ldots, n , 0 \le t \le T,}
\end{eqnarray}
respectively. Here the function $ U\dvtx[0, \infty) \times(0, \infty)^n
\rightarrow(0,1] $ is the smallest nonnegative solution, in the class
$ \mathcal{C}^{1,2} $, of the fully nonlinear parabolic partial
differential inequality
%
\begin{equation}
\label{E.3.in}
{ \partial U \over\partial\tau} (\tau, \mathbf{z})
\ge\widehat{\mathcal{L}} U (\tau, \mathbf{z}) ,\qquad
(\tau, \mathbf{z}) \in(0, \infty) \times(0, \infty)^n
\end{equation}
subject to the initial condition $ U (0, \cdot) \equiv1 $, with
%
\begin{equation}
\label{E.3.d}
\widehat{\mathcal{L}}f (\mathbf{ z}) = \sup_{a \in\mathcal{A}
(\mathbf
{ z})} \sum_{i=1}^n
\sum_{j=1}^n z_i z_j a_{ij} \biggl( {1 \over2 } D^2_{ij}f (\mathbf{ z}) + {
D_{i}f (\mathbf{ z}) \over z_1 + \cdots+ z_n } \biggr).
\end{equation}
We use in (\ref{FG}), (\ref{E.3.d}) and throughout this paper, the
notation $ D_i f = \partial f /\partial x_i $, $D^2_{ij} f = \partial^{
2} f / \partial x_i \,\partial x_j $, $ D f = (D_1 f, \ldots, D_n f )' $,
$ D^2 f = \{ D^2_{ij} f \}_{1 \le i, j \le n} $ and define
%
\begin{equation}
\label{E.3.a}
\mathcal{A} (\mathbf{y}) := \{ a \in\mathbb{S}^n \dvtx
(\theta, a) \in\mathcal{K}(\mathbf{y}) , \mbox{ for some } \theta
\in\mathbb{R}^n \} ,\qquad  \mathbf{ y} \in\mathfrak{S}_n .
\end{equation}

We call the function $ U (\cdot, \cdot) $ the \textit{arbitrage
function,} as $ U (T, \mathbf{x} ) (x_1 + \cdots+ x_n ) $ gives the
smallest initial capital starting with which an investor, who uses
nonanticipative investment rules, can match or outperform the market
portfolio by time $ t=T $, if the initial configuration of asset
capitalizations is $ \mathfrak{X} (0)= \mathbf{ x}= (x_1 , \ldots,
x_n )'
\in(0, \infty)^n $ at $ t=0 $, and does so with probability one under
any ``admissible'' market configuration that might materialize. It is
perhaps worth noting that this function $ U (\cdot, \cdot) $ is
characterized almost entirely in terms of the prevalent covariance
structure $ \alpha$. The relative risk $ \vartheta$ enters only
indirectly, namely, in determining the family of sets (\ref{E.3.a})
which are admissible for the covariance structure. Put a bit
differently, the only role~$ \vartheta$ plays is to ensure the asset
capitalization process $ \mathfrak{X} (\cdot) = ( X_1 (\cdot),
\ldots,
X_n (\cdot) )' $ takes values in $ (0, \infty)^n $.

Under additional regularity conditions, $ U (T, \mathbf{ x} ) $ is the
value of a stochastic control problem: the maximal probability that the
diffusion process $
\mathfrak{Y}(\cdot)= ( Y_1(\cdot), \ldots, Y_n (\cdot) )' $
with initial configuration $ \mathfrak{Y} (0) =\mathfrak{X} (0)=
\mathbf
{ x} \in(0, \infty)^n $, values in the punctured nonnegative orthant $
\mathfrak{S}_n $ of (\ref{1.b}), infinitesimal generator
\[
\sum_{i=1}^n \sum_{j=1}^n y_i y_j \mathrm{a}_{ij} (t, \mathbf{y}) \biggl( {1
\over2 }D^2_{ij}f (\mathbf{y}) + {D_{i}f (\mathbf{y}) \over y_1 +
\cdots+ y_n } \biggr)
\]
and controlled through the choice of covariance function $ \mathrm{a}
\dvtx
[0, \infty) \times\mathfrak{S}_n \rightarrow\mathbb{S}^n $ which
satisfies $ \mathrm{a} (t, \mathbf{y}) \in\mathcal{A} (\mathbf{y}) $
for all $ (t, \mathbf{y}) \in[0, \infty) \times\mathfrak{S}_n $,
does not hit the boundary of the orthant $ [0,\infty)^n $ by time $ t =
T $. Under appropriate conditions the function $ U (\cdot, \cdot) $
satisfies then, in the notation of (\ref{E.3.d}), the \textit{Hamilton--Jacobi--Bellman} (\textit{HJB}) \textit{equation}
%
\begin{equation}
\label{E.3.c}
( \partial/ \partial\tau) U (\tau, \mathbf{z})
= \widehat{\mathcal{L}} U (\tau, \mathbf{z}) \qquad  \mbox{on } (0,
\infty)
\times(0, \infty)^n .
\end{equation}
\textit{Relation to extant work}:
Stochastic control problems of
the ``maximal probability of containment'' type were apparently
pioneered by Van Mellaert and Dorato (\citeyear{VanDOR72}) (see also \citet{FleRis75}, pages 157--158). The ``Knightian uncertainty''
constraint imposed in (\ref{E.3}) is very similar to the formulation of
stochastic control and stochastic game problems for one-dimensional
diffusions pioneered by William Sudderth, that Sudderth and his
collaborators developed in a series of articles that
includes \citet{PesSud85}, Heath et al. (\citeyear{Heaetal87}), \citet{OrePesSud87}, \citet{SudWee89}; indeed, the developments
in Sections \ref{sec6}--\ref{sec8} of our paper can be construed as a multidimensional
extension of the Sudderth approach.

We rely strongly on Krylov's (\citeyear{Kry89N1}, \citeyear{Kry02}) work, which studies solutions
of stochastic differential equations with constraints on the drift and
diffusion coefficients in terms of ``supermartingale problems'' and
characterizes sets of stochastic integrals via appropriate supermartingales.

The approach we adopt has a lot in common with the effort, started in
the mid-1990s, to understand option pricing and hedging in the presence
of uncertainty about the underlying volatility structure of assets. We
have been influenced by this strand of work, particularly by the papers
of \citet{LYO95}, \citet{RomVar00}, \citet{GozVar02}, \citet{VAR}, \citet{TalZhe02}; other important papers
include \citet{AVELEVPAR95}, \citet{ElKJeaShr98},
\citet{CviPhaTou99}, \citet{Fre00}, Ekstr\"{o}m and Tysk (\citeyear{EksTys04}),
\citet{Mey06}, \citet{DenMar06}, whereas the recent preprints by
\citeauthor{SONTOUZHAN1} (\citeyear{SONTOUZHAN1}, \citeyear{SONTOUZHAN2}) contain very relevant results. Similar in this
spirit is the strand of work by Shige Peng and his collaborators,
surveyed in \citet{PEN}, regarding the so-called ``$ G$-Brownian
motion'' which exhibits volatility uncertainty [see also \citet{VOR}, as well as \citet{NUT} for extensions to settings where the
range of uncertainty is stochastic]. Whereas in both these strands the
relevant fully nonlinear parabolic-type partial differential (so-called
``Black--Scholes--Barenblatt'')\break equation has a typically unique
solution, here the main interest arises from \textit{lack of uniqueness}
on the part of the rather similar, fully nonlinear equation~(\ref{E.3.c}).

Let us mention that optimization problems in stochastic control,
mathematical economics and finance that involve model uncertainty have
also been treated by other authors, among them \citet{GilSch89},
\citet{GUN05}, Shied and Wu (\citeyear{SchWu05}),
\citet{KarZam05}, F\"{o}llmer and Gundel (\citeyear{FolGun06}),
\citet{Sch07}, \citet{Rie09}, Bayraktar and Yao (\citeyear{BAYYAO}),
Bayraktar, Karatzas and Yao (\citeyear{BAYKARYAO}) and \citet{KARROB} [see also the survey by
\citet{FOLSCHWEB}].

\textit{Preview}: Sections \ref{sec2} and \ref{sec3} set
up the
model for an equity market with Knightian model uncertainty regarding
its volatility and market-price-of-risk characteristics, and for
investment rules in its context. Section \ref{sec4} introduces the
notion of optimal arbitrage in this context, whereas Section \ref{sec5}
discusses the relevance of the fully nonlinear parabolic partial
differential inequality of HJB type (\ref{E.3.in}), (\ref{E.3.d}) in
characterizing the arbitrage function and in finding an investment rule
that realizes the best outperformance of the market portfolio. Section
\ref{sec6} presents a verification-type result for this equation.
Sections \ref{sec7} and \ref{sec8} make the connection with the
stochastic control problem of maximizing the probability of containment
for an auxiliary It\^{o} process, controlled in a nonanticipative way
and in a manner that affects both its drift and dispersion
characteristics. Finally, Section \ref{sec9} develops yet another
characterization of the arbitrage function, this time as the min-max
value of a zero-sum stochastic game; the investment rule that realizes
the best outperformance of (optimal arbitrage with respect to) the
market, is now seen as the investor's best response to a ``least
favorable'' market model configuration.


\section{Equity market with Knightian model uncertainty}
\label{sec2}

We shall fix throughout a canonical, filtered measurable space $
(\Omega, \mathcal{F} )$, $ \mathbb{F} = \{ \mathcal{F}(t)\}_{0
\le t <\infty} $ and assume that $ \Omega$ contains the space $
\mathfrak{W} \equiv C ( [0,\infty); (0, \infty)^n ) $ of all continuous
functions $
\mathfrak{w} \dvtx[0, \infty) \rightarrow( 0,\infty)^n $. We
shall specify this canonical space in more detail in Section \ref{sec7} below,
when such detail becomes necessary.

On this space, we shall consider a vector of continuous, adapted
processes $ \mathfrak{X}(\cdot) = ( X_1(\cdot), \ldots, X_n(\cdot)
)' $ with values in $ (0, \infty)^n $; its components will represent
stock capitalizations in an equity market with $ n $ assets, and thus
the total market capitalization will be the sum
%
\begin{equation}
\label{1.a}
X(t) := X_1 (t)+ \cdots+X_n (t) ,\qquad  0 \le t
< \infty.
\end{equation}
We shall also fix throughout a collection $
\mathbb{K} = \{ \mathcal{K}
(\mathbf{y}) \}_{\mathbf{y} \in\mathfrak{S}_n}
$ of nonempty, compact and convex subsets of $ \mathbb{R}^n \times
\mathbb{S}^n $ as in (\ref{1.b}).

We shall consider $ \mathbb{R}^n$-valued functionals $ \vartheta
(\cdot
, \cdot) = ( \vartheta_1(\cdot, \cdot),
\ldots, \vartheta_n (\cdot, \cdot) )'$ and $ \mathbb{S}^n$-valued
functionals $ \alpha(\cdot, \cdot) =(
\alpha_{i j}(\cdot, \cdot) )_{1 \le i, j \le n} $, all of them
defined on $ [0, \infty) \times\Omega$ and progressively measurable
[see \citet{KarShr91}, Definition~3.5.15]. We shall assume
that, for every continuous function $ \mathfrak{w}\dvtx[0, \infty)
\rightarrow( 0,\infty)^n $ and $ T \in(0,\infty) $, these functionals
satisfy the constraint and integrability conditions, respectively,
\begin{eqnarray}
\label{1.2}
( \vartheta(T, \mathfrak{w} ), \alpha(T, \mathfrak{w}) ) &\in&
\mathcal
{K} ( \mathfrak{w} (T) ) , \nonumber\\[-8pt]\\[-8pt]
\int_0^T \bigl( \| \vartheta(t, \mathfrak{w})
\|^2 + \operatorname{Tr} ( \alpha(t, \mathfrak{w}) ) \bigr) \,\mathrm{d} t &<&
\infty.\nonumber
\end{eqnarray}
We shall also consider $(n \times
n)$-matrix-valued functionals
\[
\sigma(\cdot, \cdot) = (
\sigma_{i \nu}(\cdot, \cdot) )_{1 \le i, \nu\le n} ,
\]
where
%
\begin{equation}
\label{AlphaSigma}
\qquad \sigma(t , \mathfrak{w}) = \sqrt{\alpha(t , \mathfrak{w}) } \mbox{ is
a square root of } \alpha( t , \mathfrak{w})\dvtx
\alpha( t , \mathfrak{w}) = \sigma
( t , \mathfrak{w}) \sigma' ( t , \mathfrak{w}) .\hspace*{-10pt}
\end{equation}


\subsection{Admissible systems}\label{subsec2.1}

For a given collection of sets $ \mathbb{K} $ as in (\ref{1.b})
and a~fixed initial configuration $ \mathbf{x} = (x_1, \ldots,
x_n)' \in(0, \infty)^n $ of asset capitalizations, we shall call
\textit{admissible system} a collection $ \mathcal{M} $ consisting of the
underlying filtered space $ (\Omega, \mathcal{F})$, $ \mathbb{F} = \{
\mathcal{F}(t)\}_{0 \le t <\infty} $, of a probability measure $
\mathbb
{P} $ on it, and of a pair of processes $ ( \mathfrak{X}(\cdot),
W(\cdot
) ) $, with $ W(\cdot) = ( W_1(\cdot), \ldots, W_n (\cdot) )' $ an $
n$-dimensional $ \mathbb{F}$-Brownian motion
and $ \mathfrak{X}(\cdot)= ( X_1(\cdot), \ldots, X_n (\cdot) )' $
taking values in $ (0, \infty)^n $. These processes have the dynamics
%
\begin{equation}
\label{1.1}
\qquad \mathrm{d} X_i(t) = X_i(t)
\sum_{\nu=1}^n \sigma_{i \nu}(t, \mathfrak{ X} ) [ \mathrm{d} W_\nu(t) +
\vartheta_\nu(t, \mathfrak{ X})\, \mathrm{d} t ] ,\qquad  X_i(0)=x_i>0
\end{equation}
for some progressively measurable functionals $ \vartheta(\cdot,
\cdot
) $ and $ \sigma(\cdot, \cdot) $ which satisfy~(\ref{1.2}) and (\ref
{AlphaSigma}) above.

We shall think of this admissible system $ \mathcal{M} $ as a \textit
{model} subject to ``Knightian'' uncertainty; this is
expressed by the requirement $ ( \vartheta( t, \mathfrak{
X}), \alpha( t, \mathfrak{ X}) ) \in\mathcal{K} (
\mathfrak{X} ( t)
) $ in (\ref{1.2}), (\ref{1.b}) about the market
price of risk and the covariance structure of the asset capitalization
vector process $ \mathfrak{X}(\cdot) $. In order not to lose sight of
the underlying
probability space, we shall denote by
$ \mathbb{P}^\mathcal{M} $ (resp., $
\mathbb{E}^{ \mathbb{P}^\mathcal{M}} $) the probability measure (resp.,
the corresponding expectation operator) on this space. Finally, $
\mathfrak{M} (\mathbf{x}) $ will denote the collection of all such
admissible systems or ``models'' with initial configuration $ \mathbf{x}
= (x_1, \ldots, x_n)' \in(0, \infty)^n $. We shall think of $
\mathfrak
{M} (\mathbf{x}) $ as a \textit{meta-model}, a collection of admissible
models, and of the collection $ \mathbf{M}= \{ \mathfrak{M} (\mathbf
{x}) \}_{\mathbf{x}\in(0, \infty)^n} $ as a ``family of meta-models.''

The interpretation is that the components of the
driving Brownian motion~$ W(\cdot) $ represent the independent
factors of the resulting model; the entries of the matrix $
\sigma(t, \mathfrak{ X}) $ are the local volatility rates of the
asset capitalization vector process $ \mathfrak{X}(\cdot) $ at time $
t $; the entries of the matrix $ \alpha(t, \mathfrak{ X}) $ as in~%
(\ref
{AlphaSigma}) represent the local covariance rates; whereas the
components of the vector~$ \vartheta(t, \mathfrak{ X}) $ are the
market price of risk (also called \textit{relative risk}) rates prevalent
at time
$ t $. In particular,
%
\begin{equation}
\label{b}
\beta(t, \mathfrak{ w}) := \sigma(t, \mathfrak{w}) \vartheta(t,
\mathfrak{w}) = \sqrt{\alpha(t, \mathfrak{w}) } \vartheta(t,
\mathfrak{w})
\end{equation}
is, in the notation of (\ref{AlphaSigma}), the vector of mean rates of
return for the various assets at time $ t $, when the equations of
(\ref
{1.1}) are cast in the more familiar form
%
\begin{equation}
\label{b1}
\mathrm{d} X_i(t) = X_i(t) \Biggl( \beta_i (t, \mathfrak{ X}) \,\mathrm{d}
t +
\sum_{\nu=1}^n \sigma_{i \nu}(t, \mathfrak{ X}) \,\mathrm{d} W_\nu
(t) \Biggr)
,\quad  i=1, \ldots, n .\hspace*{-35pt}
\end{equation}

The integrability condition of (\ref{1.2}) guarantees that the process
$ \mathfrak{X}(\cdot) $ takes values in $ (0, \infty)^n $, $ \mathbb
{P}$-a.s.; it implies also that the exponential process
%
\begin{equation}
\label{expo}
L (t) := \exp\Biggl\{{-}\int_0^t
\vartheta' (s, \mathfrak{ X}) \,\mathrm{d} W(s) - { 1 \over2} \int_0^t
\| \vartheta(s, \mathfrak{ X}) \|^2 \,\mathrm{d} s \Biggr\} ,\quad  0 \le t <
\infty\hspace*{-35pt}
\end{equation}
is well defined and a strictly positive local martingale, thus also a
supermartingale.

This process $ L(\cdot) $ plays the role of a state-price-density or
``deflator'' in the present context. Just as in our earlier works,
Fernholz and Karatzas (\citeyear{FerKar10}, \citeyear{FerKar10b}) as well as \citet{RUF}---mostly
in a Markovian context, and without model uncertainty---an important
feature of this subject is that $ L(\cdot) $ has to be allowed to be a
\textit{strict} local martingale, that is, that $ \mathbb{E}^\mathbb{P}
(L(T))<1 $ be allowed to hold for some, if not all, $ T \in(0, \infty) $.

\subsection{Supermartingale problems}
\label{subsec2.2}

Constraint (\ref{E.3}) brings us in the realm of the Krylov (\citeyear{Kry89N1}, \citeyear{Kry02}) approach, which studies stochastic differential equations with
constraints on the drift and diffusion coefficients in terms of
``supermartingale problems.'' In particular, Theorem 2.2 of \citet{Kry02} shows that, under a suitable regularity condition on the family
of sets~$ \mathbb{K} $ in~(\ref{1.b}), solving stochastic equation~%
(\ref
{1.1}) subject to the requirements of~(\ref{1.2}) can be cast as a
\textit{supermartingale problem,} as follows.

Consider the nonlinear partial differential operator associated with
(\ref{b1}), (\ref{b}), namely
%
\begin{eqnarray}
\label{LF}
\mathcal{L}f (\mathbf{ z}) &\hspace*{3pt}=& F ( D^2 f (\mathbf{ z}), D f (\mathbf{
z}), \mathbf{ z} )  \nonumber
\\
\eqntext{\mbox{with }\displaystyle F (Q, p, \mathbf{z}) :=
\mathop{\sup_{(\theta, a) \in\mathcal{K} (\mathbf{ z}) }}_{b = \sqrt{a }
\theta} \Biggl( { 1 \over2 } \sum_{i=1}^n
\sum_{j=1}^n z_i z_j a_{ij} Q_{ij} + \sum_{i=1}^n z_i b_i p_{i} \Biggr) ,} \\[-38pt]
\\[10pt]
\eqntext{(Q,p, \mathbf{ z}) \in \mathbb{S}^n \times\mathbb{R}^n \times(0,
\infty
)^n .}
\end{eqnarray}
The supermartingale problem is to find a probability measure $ \mathbb
{P} = \mathbb{P}^\mathcal{M} $ on the filtered measurable space $
(\Omega, \mathcal{F} )$, $ \mathbb{F} = \{ \mathcal{F}(t)\}_{0
\le t <\infty} $, under which $ \mathfrak{X}(\cdot) $
takes values in $ (0, \infty)^n $ a.s., and the process
\[
u (t , \mathfrak{X} (t) ) - \int_0^t \biggl( { \partial u \over\partial s }
(s , \mathfrak{X} (s) ) + \mathcal{L} u ( \mathfrak{X} (s) ) \biggr)
\,\mathrm{d}s ,\qquad  0 < t < \infty
\]
is a local supermartingale for every $ u \dvtx(0, \infty) \times(0,
\infty
)^n \rightarrow\mathbb{R} $ of class $ \mathcal{C}^{1,2} $ with
compact support.

The regularity condition on the family $ \mathbb{K} $ in (\ref{1.b})
that we alluded to earlier, mandates that the function
%
\begin{equation}
\label{measur}
F (Q,p, \cdot) \mbox{ in (\ref{LF}) is Borel measurable, for every }
(Q, p) \in\mathbb{S}^n \times\mathbb{R}^n .\hspace*{-25pt}
\end{equation}
If, in addition, the family of sets $ \mathbb{K} $ in (\ref{1.b})
satisfies the linear growth condition
%
\begin{equation}\label{LinGro}
\mathop{\sup_{(\theta, a) \in\mathcal{K} (\mathbf{y}) }}_{ b = \sqrt{a }
\theta} \Biggl( \sum_{i=1}^n \sum_{i=1}^n
y_i y_j a_{ij} + \sum_{i=1}^n (y_i b_i)^2 \Biggr)^{1/2}
\le C ( 1 + \| \mathbf{ y} \| )
\end{equation}
and the upper-semicontinuity condition
%
\begin{equation}
\label{UppSemiCont}
\limsup_{ [0,\infty)^n
\ni\mathbf{z} \rightarrow\mathbf{ y} } F (Q, p , \mathbf{z}) \le
F (Q, p, \mathbf{ y}) \qquad  \forall(Q, p) \in\mathbb{S}^n \times
\mathbb{R}^n
\end{equation}
for every $ \mathbf{ y} \in\mathfrak{S}_n $ and some real constant $
C>0 $, then Theorem 3.2 in \citet{Kry02} shows that the family $ \{
\mathbb{P}^\mathcal{M} \}_{\mathcal{M} \in\mathfrak{ M}(\mathbf{ x})}
$ is convex and sequentially compact in the topology of vague
convergence of probability measures.


\subsection{Markovian admissible systems}
\label{subsec2.3}
We shall also consider the subcollection $ \mathfrak{M}_* (\mathbf{x})
\subset\mathfrak{M} (\mathbf{x}) $ of \textit{Markovian} admissible
systems, for which the functionals $ \vartheta(\cdot, \cdot) $ and $
\alpha(\cdot, \cdot) $ as in (\ref{1.2})--(\ref{1.1}) are given as
%
\begin{equation}
\label{marko}
\vartheta(t, \mathfrak{X}) = \theb(t, \mathfrak{X} (t)) ,\qquad  \alpha
(t, \mathfrak{X}) = \mathbf{a} ( t, \mathfrak{X} (t)) ,
\end{equation}
with measurable functions $ \theb\dvtx[0, \infty) \times
(0, \infty)^n \rightarrow\mathbb{R}^n $, \mbox{$ \mathbf{a} \dvtx[0,
\infty)\times(0, \infty)^n \rightarrow\mathbb{S}^n $} that satisfy
%
\begin{equation}
\label{D.22}
( \theb( t, \mathbf{ z}), \ab( t, \mathbf{ z})
) \in\mathcal{K} (\mathbf{ z})  \qquad \forall
(t, \mathbf{ z}) \in[0, \infty) \times
(0, \infty)^n .
\end{equation}
%
Under condition (\ref{LinGro}) it follows then from so-called
Markovian selection results [\citet{Kry73}; \citet{StrVar79}, Chapter 12; \citet{EthKur86}, Section 4.5] that the
state process $ \mathfrak{ X} (\cdot) $ of (\ref{1.1}) can be assumed
to be (strongly) Markovian under $ \mathbb{P}^{\mathcal{M}} $, $
\mathcal{M} \in\mathfrak{M}_* (\mathbf{x}) $. We shall make this
selection whenever admissible systems in $ \mathfrak{M}_* (\mathbf{x})
$ are invoked.


\section{Investment rules}
\label{sec3}

Consider now an investor who is ``small'' in the sense that his actions
do not affect market prices. He starts with initial fortune $ v>0 $ and uses
a rule that invests a proportion $ \varpi_i (t) = \Pi_i ( t,
\mathfrak{X} ) $ of current wealth in the $ i$th
asset of the equity market, for any given time $ t \in[0 , \infty) $
and all $ i=1, \ldots, n $; the
remaining proportion $ \varpi_0 (t) := 1 - \sum_{i=1}^n \varpi_i
(t) $ is held in cash (equivalently, in a zero-interest
money market). Here $ \Pi\dvtx[0, \infty) \times\mathfrak{W}
\rightarrow\mathbb{R}^n $ is a progressively measurable functional
assumed to satisfy the requirement
%
\begin{eqnarray}
\label{1.3}
\int_0^T \bigl( | \Pi^{ \prime} ( t, \mathfrak{w}
) \sigma(t, \mathfrak{w}) \vartheta(t, \mathfrak{w}) | + \Pi^{
\prime} ( t, \mathfrak{w} ) \alpha(t, \mathfrak{w}) \Pi( t, \mathfrak{w} )
\bigr) \,\mathrm{d} t < \infty,\nonumber\\[-8pt]\\[-8pt]
\eqntext{\forall T \in(0, \infty)}
\end{eqnarray}
for every continuous function $ \mathfrak{w} \dvtx[0, \infty)
\rightarrow (0, \infty)^n $. (Thus, the requirement (\ref{1.3}) will be
in force under all admissible systems.) We shall denote throughout by $
\mathfrak{P} $ the collection of all such (nonanticipative)
\textit{investment rules}, and by $ \mathfrak{P}_* $ the sub-collection
of all \textit{Markovian} investment rules, that is, those that can be
expressed as $ \varpi_i (t) = \pi_i ( t, \mathfrak {X}( t) ) $, $ 0 \le
t < \infty, i=1, \ldots, n $ for some measurable function $ \pi\dvtx[0,
\infty) \times(0, \infty)^n\rightarrow\mathbb{R}^n $.

An investment rule is called \textit{bounded}, if the functional $ \Pi$
is bounded uniformly on $ [0,\infty) \times\mathfrak{W} $; for a
bounded investment rule, requirement
(\ref{1.3}) is satisfied automatically, thanks to (\ref{1.2}). An
investment rule is called \textit{portfolio} if the functional $
\Pi$ satisfies $ \sum_{i=1}^n \Pi_i =1 $ on $ [0,\infty)
\times\mathfrak{W} $, and a portfolio is called \textit{long-only}
if $ \Pi_1 \ge0, \ldots, \Pi_n \ge0 $ also hold
on this domain. A long-only portfolio is clearly bounded.

Given an initial wealth $ v \in(0, \infty) $, an investment
rule $ \Pi\in\mathfrak{ P} $ and an admissible model $
\mathcal{M} \in\mathfrak{M} (\mathbf{x}) $, the resulting wealth
process $ Z(\cdot) \equiv Z^{ v, \Pi} (\cdot) $ satisfies the
dynamics
%
\begin{equation}
\label{1.5}
\qquad\qquad \frac{ \mathrm{d} Z (t)}{Z (t)} =
\sum_{i=1}^n \Pi_i ( t, \mathfrak{X} )\,
\frac{ \mathrm{d} X_i(t)}{X_i(t)} = \Pi^{ \prime}
( t, \mathfrak{X}
) \sigma(t, \mathfrak{X})
[ \vartheta(t, \mathfrak{X}) \,\mathrm{d} t + \mathrm{d} W (t)
]
\end{equation}
and the initial condition $ Z(0) =v $. In conjunction with
(\ref{expo}) in the differential form
%
\begin{equation}
\label{expotoo}
\mathrm{d} L ( t) = - L ( t) ( \vartheta(t, \mathfrak{ X}) )'\,\mathrm{d}
W(t)\vadjust{\goodbreak}
\end{equation}
and the product rule of the stochastic calculus, this
gives
%
\begin{eqnarray}
\label{1.6}
L ( t) Z^{ v, \Pi} ( t) = v + \int_0^t L (s) Z^{ v, \Pi} (s) \bigl(
\sigma' (s, \mathfrak{X}) \Pi( s, \mathfrak{X} ) - \vartheta(s,
\mathfrak{X})
\bigr)' \,\mathrm{d} W (s) ,\nonumber\\[-8pt]\\[-8pt]
\eqntext{0 \le t< \infty.}
\end{eqnarray}
For any initial configuration $ \mathbf{x} =
(x_1, \ldots, x_n)' \in(0, \infty)^n $, initial wealth $ v
\in(0, \infty) $, investment rule $ \Pi\in\mathfrak{P} $ and
admissible model $ \mathcal{M} \in\mathfrak{M} (\mathbf{x}) $,
the product $ L( \cdot) Z^{ v, \Pi} ( \cdot) $ is therefore under $
\mathbf{ }
\mathbb{P}^\mathcal{M} $ a continuous,
positive local martingale, thus also a supermartingale. Once again, it
is important that this process be allowed to be a strict local martingale.

\subsection{The market portfolio}
\label{subsec3.1}

The choice of Markovian investment rule $ \mathfrak{m} \in\mathfrak{
P}_* $ given by
\[
\mathfrak{m}_i (t,\mathbf{z}) = { z_i \over z_1 + \cdots+ z_n } ,\qquad
i=1, \ldots,n ,\ t \in[0, \infty) ,\ \mathbf{z} \in(0,\infty)^n
\]
leads to the long-only \textit{market portfolio}
%
\begin{eqnarray}
\label{1.7.a}
\mu(t) = \mathfrak{ X} (t) / X(t) \nonumber\\[-8pt]\\[-8pt]
\eqntext{\mbox{with weights }
\mu_i (t) = X_i (t) / X(t) ,\ i=1, \ldots, n ,\ 0 \le t
< \infty}
\end{eqnarray}
in the notation of (\ref{1.a}). It follows from (\ref{1.5}) that
investing according
to this portfolio amounts to owning the entire market, in
proportion to the initial wealth, $ Z^{v, \mathfrak{m}} (\cdot) = v
X(\cdot) / X(0) $.

\subsection{Ramifications}
\label{subsec3.2}

Reading (\ref{1.6}) for the market portfolio of (\ref{1.7.a}), and
recalling (\ref{1.a}), leads to
%
\begin{eqnarray}
\label{1.c}
L ( t) X ( t) = X(0) + \int_0^t L (s) \bigl(
\sigma' (s, \mathfrak{ X}) \mathfrak{X} (s) - \vartheta(s,
\mathfrak{
X}) X(s) \bigr)' \,\mathrm{d} W (s) ,\nonumber\\[-8pt]\\[-8pt]
\eqntext{0 \le t< \infty}
\end{eqnarray}
or equivalently $ \mathrm{d} ( L ( t) X ( t) ) = - L ( t) X ( t) (
\widetilde{\vartheta} (t, \mathfrak{ X}) )'\, \mathrm{d}W(t) $, where
%
\begin{eqnarray}
\label{1.d}
\widetilde{\vartheta} (t, \mathfrak{w}) := \vartheta(t, \mathfrak{w})
- { \sigma' (t, \mathfrak{w}) \mathfrak{w} (t) \over\mathfrak{w}_1
(t) + \cdots+ \mathfrak{w}_n (t) } \nonumber\\[-8pt]\\[-8pt]
\eqntext{\displaystyle \mbox{satisfies }
\int_0^T \| \widetilde{\vartheta} (t, \mathfrak{w}) \|^2 \,\mathrm{d} t <
\infty}
\end{eqnarray}
for all $ (t, \mathfrak{w} ) \in[0, \infty) \times\mathfrak{W} $,
thanks to (\ref{1.2}). With this notation, it follows from (\ref
{1.c}) that
\[
L (\cdot) X (\cdot) = (x_1 + \cdots+ x_n) \cdot
\exp\Biggl\{ - \int_0^\cdot( \widetilde{\vartheta} (t, \mathfrak{ X} )
)' \,\mathrm{d} W(t) - { 1 \over2} \int_0^\cdot\| \widetilde
{\vartheta
} (t, \mathfrak{ X} ) \|^2\, \mathrm{d} t \Biggr\} .
\]
On the strength of the integrability condition in (\ref{1.d}), the
Dambis--Dubins--Schwartz representation [e.g., \citet{KarShr91}, page 174] of the $ \mathbb{P}^\mathcal{M}$-local martingale
%
\begin{eqnarray}
N (\cdot) := \int_0^\cdot( \widetilde{\vartheta} (t, \mathfrak{
X} ))'\, \mathrm{d} W(t)\nonumber \\
\eqntext{\displaystyle \mbox{with quadratic variation } \langle N \rangle(
\cdot) = \int_0^\cdot\| \widetilde{\vartheta} (t, \mathfrak{ X} )
\|^2 \,\mathrm{d} t < \infty}
\end{eqnarray}
gives
%
\begin{equation}
\label{repr}
\quad\qquad L (T) X (T) = (x_1 + \cdots+ x_n) \cdot e^{ B(u) - (u/2)} |_{ u =
\langle N \rangle(T)} ,\qquad  0 \le T < \infty,
\end{equation}
where $ B(\cdot) $ is one-dimensional, standard Brownian motion under $
\mathbb{P}^\mathcal{M} $.\break Whereas
the equations of (\ref{1.1}) can be written as
%
\begin{equation}
\label{1.1.a}
\quad\mathrm{d} X_i(t) = X_i(t)
\Biggl( { \sum_{j=1}^n \alpha_{ij} (t, \mathfrak{ X} ) X_j (t) \over X_1
(t) + \cdots+ X_n (t) } \,\mathrm{d} t + \sum_{\nu=1}^n \sigma_{i
\nu
}(t, \mathfrak{ X}) \,\mathrm{d} \widetilde{W}_\nu(t) \Biggr)
\end{equation}
for $ i=1, \ldots, n $, with
%
\begin{equation}
\label{1.1.b}
\widetilde{W} (\cdot) := W(\cdot) + \int_0^\cdot\widetilde
{\vartheta
} ( t, \mathfrak{ X} ) \,\mathrm{d} t .
\end{equation}
We then have the representation
\begin{eqnarray}\label{Lambda}
\qquad\Lambda(\cdot) &:=& { X (0) \over L (\cdot) X (\cdot) }
=
\exp\Biggl\{ \int_0^\cdot( \widetilde{\vartheta} (t, \mathfrak{ X} ) )'
\,\mathrm{d} W(t) + { 1 \over2} \int_0^\cdot\| \widetilde{\vartheta}
(t, \mathfrak{ X} ) \|^2\, \mathrm{d} t \Biggr\}
\nonumber\\[-8pt]\\[-8pt]
&\qquad\hspace*{3pt}=&
\exp\Biggl\{ \int_0^\cdot
( \widetilde{\vartheta} (t, \mathfrak{ X} ) )' \,\mathrm{d}
\widetilde
{W}(t) - { 1 \over2} \int_0^\cdot
\| \widetilde{\vartheta} (t, \mathfrak{ X} ) \|^2 \,\mathrm{d} t \Biggr\}\nonumber
\end{eqnarray}
for the normalized reciprocal of the deflated total market capitalization,
and
%
\begin{equation}
\label{MU}
\mathrm{d} \mu_i (t) = \mu_i (t) \bigl( \mathfrak{e}_i - \mu(t) \bigr)'
\sigma
(t, \mathfrak{ X} ) \,\mathrm{d} \widetilde{W} (t) ,\qquad  i=1, \ldots, n
\end{equation}
for the dynamics market weights in (\ref{1.7.a}); here $ \mathfrak{e}_i
$ is the $ i$th unit vector in~$ \R^n $.


\section{Optimal arbitrage relative to the market}
\label{sec4}

Let us consider now the smallest proportion
\begin{eqnarray}
\label{1.7}
\qquad\mathfrak{u} (T,\mathbf{x}) &=& \inf\bigl\{ r>0 \dvtx\exists
\Pi_r \in\mathfrak{ P}, \mbox{ s.t. } \mathbb{P}^\mathcal{M} \bigl( Z^{ r
X(0), \Pi_r} ( T) \ge X(T) \bigr) = 1 ,\nonumber\\[-8pt]\\[-8pt]
&&\qquad\hspace*{198pt}\forall\mathcal{M} \in\mathfrak{M} (\mathbf{x}) \bigr\}\nonumber
\end{eqnarray}
of the initial total market capitalization $ X(0) = x_1 + \cdots
+x_n $ which allows the small investor, starting with initial capital $
\mathfrak{u} (T,\mathbf{x}) X(0) $ and through judicious
choice of investment rule in the class $ \mathfrak{ P} $, to match or
exceed the performance of the market portfolio over the time-horizon $
[0,T] $, and to do this with $ \mathbb{P}^\mathcal{M}$-probability one,
under any model $ \mathcal{M} \in\mathfrak{M} (\mathbf{ x}) $ that
might materialize. We shall refer to $ \mathfrak{u} (\cdot, \cdot) $
of (\ref{1.7}) as the \textit{arbitrage function} for the family of
meta-models $ \mathbf{M} = \{ \mathfrak{M} (\mathbf{x}) \}_{\mathbf
{x}\in(0, \infty)^n} $, and think of it as a version of the arbitrage
function studied in Fernholz and Karatzas (\citeyear{FerKar10}) which is ``robust''
with respect to $ \mathbf{M} $.

The quantity of (\ref{1.7}) is strictly positive; see Proposition \ref{Proposition_1}
below and the discussion following it. On the other hand, the set of
(\ref{1.7}) contains the number $ r=1 $, so clearly
\[
0 < \mathfrak{u} (T,\mathbf{x}) \le1 .
\]

If $ \mathfrak{u} (T,\mathbf{x}) <1 $, then for every $ r \in(
\mathfrak{u} (T,\mathbf{x}), 1 ) $---and even for $ r = \mathfrak{u}
(T,\mathbf{x}) $ when the infimum in (\ref{1.7}) is attained, as indeed
it is in the context of Theorem 1 below---there exists an investment
rule $ \Pi_r \in\mathfrak{ P} $ such that
\[
Z^{ X(0), \Pi_r} (T) \ge{1 \over r } X(T) > X(T) = Z^{ X(0), \mathrm
{m}} (T) ,\qquad  \mathbb{P}^\mathcal{M}\mbox{-a.s.}
\]
holds for every $ \mathcal{M} \in\mathfrak{ M} (\mathbf{ x}) $. In
other words, the investment rule $ \Pi_r $ leads then to \textit{strong
arbitrage relative to the market} portfolio in the terminology of
\citet{FERKAR}---here with the extra feature that such
arbitrage is now \textit{robust}, that is, holds under any possible
admissible system or ``model'' that might materialize. If, on the other
hand, $ \mathfrak{u} (T,\mathbf{x}) =1 $, then such outperformance of
(equivalently, strong arbitrage relative to) the market is just not
possible over all meta-models $ \mathcal{M} \in\mathfrak{M} (\mathbf
{x}) $. In either case, the highest return on investment relative to
the market
\begin{eqnarray*}
\mathfrak{b} (T,\mathbf{x}) &:=& \sup\bigl\{ b>0 \dvtx\exists
\Pi\in\mathfrak{ P}, \mbox{ s.t. } \mathbb{P}^\mathcal{M} \bigl( Z^{ X(0),
\Pi} ( T) \ge bX(T) \bigr) = 1 ,\\
&&\hspace*{195pt}\forall\mathcal{M} \in\mathfrak{M} (\mathbf{x}) \bigr\} ,
\end{eqnarray*}
achievable using (nonanticipative) investment rules, is given as $
\mathfrak{b} (T,\mathbf{x})= 1 / \mathfrak{u} (T,\mathbf{x}) \ge1 $.

\begin{remark}
\label{Remark_1}
Instances of $ \mathfrak{u} (T,\mathbf{x}) <1 $ occur, when there
exists a constant $ \zeta>0 $ such that either
%
\begin{equation}
\label{3.13}
\inf_{a \in\mathcal{A} (\mathbf{z})
} \Biggl( \sum_{i=1}^n { z_i a_{ii} \over z_1 + \cdots+ z_n } - \sum
_{i=1}^n \sum_{j=1}^n { z_i z_j a_{ij} \over(z_1 + \cdots+ y_n )^2 }
\Biggr) \ge\zeta
\end{equation}
or
%
\begin{equation}
\label{3.14}
\biggl( { (z_1 \cdots z_n )^{1/n} \over z_1 + \cdots+ z_n } \biggr) \cdot\inf
_{a \in\mathcal{A} (\mathbf{ z})
} \Biggl( \sum_{i=1}^n a_{ii} - { 1 \over n } \sum_{i=1}^n \sum_{j=1}^n
a_{ij} \Biggr) \ge\zeta
\end{equation}
holds for every $ \mathbf{z} \in(0,\infty)^n $. See the survey paper
\citet{FERKAR}, Examples 11.1 and 11.2 [as well as
\citet{FERKAR05}, \citet{FERKARKAR05} for additional
examples].
\end{remark}

\begin{prop}
\label{Proposition_1}
The quantity of (\ref{1.7})
satisfies
%
\begin{eqnarray}\label{1.8}
\mathfrak{u} (T,\mathbf{x}) \ge\Phi(T,\mathbf{x})>0\nonumber\\[-8pt]\\[-8pt]
\eqntext{\displaystyle \mbox{where } \Phi(T,\mathbf{x}):=
\sup_{\mathcal{M} \in
\mathfrak{M} (\mathbf{x})} \biggl( \frac{ \mathbb{E}^{ \mathbb
{P}^\mathcal{M}}
[ L(T) X(T) ] }{ x_1 + \cdots+ x_n } \biggr) .}
\end{eqnarray}
Furthermore, under conditions (\ref{measur})--(\ref{UppSemiCont}), there
exists an admissible system $ \mathcal{M}_o \in\mathfrak{M} (\mathbf
{x}) $ such that
%
\begin{equation}
\label{1.88}
\Phi(T,\mathbf{x}) = \frac{ \mathbb{E}^{ \mathbb{P}^{\mathcal{M}_o}}
[ L(T) X(T) ] }{ x_1 + \cdots+ x_n } .
\end{equation}
\end{prop}

\begin{pf}
Take an arbitrary element $ r >0 $ of the set on
the right-hand side of (\ref{1.7}) and an arbitrary admissible
system $ \mathcal{M} \in\mathfrak{M} (\mathbf{x}) $. There
exists then an investment rule $ \Pi_r \in\mathfrak{P} $ with the
inequality $ Z^{ r X(0), \Pi_r} ( T) \ge X(T) $ valid $
\mathbb{P}^\mathcal{M}$-a.s. On the strength of (\ref{expo}) and
(\ref
{1.6}), the process $ L( \cdot) Z^{ r X(0), \Pi_r} ( \cdot) $ is a~$
\mathbb{P}^\mathcal{M}$-supermartingale; thus (\ref{repr}) and (\ref
{1.d}) lead to
\begin{eqnarray}\label{1.89}
r (x_1+ \cdots+x_n) &=& r X(0) \ge\mathbb{E}^{ \mathbb{P}^\mathcal{M}}
\bigl[ L (T) Z^{ r X(0), \Pi_r} (T) \bigr] \nonumber\\[-8pt]\\[-8pt]
&\ge&\mathbb{E}^{ \mathbb{P}^\mathcal{M}} [ L (T) X(T) ] > 0 .\nonumber
\end{eqnarray}
The inequality $ \mathfrak{u} (T,\mathbf{x}) \ge\Phi(T,\mathbf{x}) $
in (\ref{1.8}) follows now from the arbitrariness of
$ r >0 $ and $ \mathcal{M} \in\mathfrak{M}
(\mathbf{x}) $. The existence of an admissible system $ \mathcal{M}_o
\in\mathfrak{M} (\mathbf{x}) $ that satisfies (\ref{1.88}) follows
from Theorem 3.4 in \citet{Kry02}, in conjunction with the dynamics of
(\ref{1.1}) and (\ref{expotoo}).
\end{pf}

Although strong arbitrage relative to the market may exist within the
framework of the models $ \mathcal{M} \in\mathfrak{M} (\mathbf{x}) $
studied here (cf. Remark \ref{Remark_1}), the existence of a strictly
positive supermartingale deflator process $ L (\cdot) $ as in (\ref
{expo}) proscribes scalable arbitrage opportunities, also known as
\textit{Unbounded Profits with Bounded Risk} (\textit{UPBR}); this is reflected in the
inequality $ \mathfrak{u} (T,\mathbf{x})>0 $ of~(\ref{1.8}). We refer
the reader to \citet{DelSch95N1} for the origin of the
resulting \textit{NUPBR} concept, and to \citet{KarKar07} for
an elaboration of this point in a different context, namely, the
existence and properties of the num\'{e}raire portfolio.

Finally, let us write (\ref{1.8}) as
%
\begin{equation}
\label{1.9}
\qquad \Phi(T,\mathbf{x}) = \sup_{\mathcal{M} \in
\mathfrak{M} (\mathbf{x})} \mathfrak{u}_\mathcal{M} (T,\mathbf{x})
\qquad \mbox{where } \mathfrak{u}_\mathcal{M} (T,\mathbf{x}) :=
\frac{ \mathbb{E}^{ \mathbb{P}^\mathcal{M}}
[ L(T) X(T) ] }{ x_1 + \cdots+ x_n } .\hspace*{-10pt}
\end{equation}
We have for this quantity the interpretation
\[
\mathfrak{u}_\mathcal{M} (T,\mathbf{x}) = \inf\bigl\{ r>0 \dvtx\exists
\Pi_r \in\mathfrak{ P}, \mbox{s.t. } \mathbb{P}^\mathcal{M} \bigl( Z^{ r
X(0), \Pi_r} ( T) \ge X(T) \bigr) = 1 \bigr\}
\]
as the ``arbitrage function for the model $ \mathcal{M} \in\mathfrak
{M}(\mathbf{x}) $,'' at least when the matrix $ \sigma(t, \mathfrak{w})
$ in invertible for every $ (t,\mathfrak{w}) \in(0, \infty) \times
\mathfrak{W} $ and when $ (\mathbb{P}^\mathcal{M}, \mathbb
{F})$-martingales can be represented as stochastic integrals with
respect to the Brownian motion $ W(\cdot) $ in (\ref{1.1}).


\section{A fully nonlinear PDI}
\label{sec5}

Consider now a continuous function $ U \dvtx[0,\infty) \times(0,
\infty)^n \rightarrow(0, \infty) $ with
%
\begin{equation}
\label{2.1}
U(0, \mathbf{z}) = 1 ,\qquad \mathbf{z} \in(0, \infty)^n ,
\end{equation}
which is of class $ \mathcal{C}^{1,2} $ on $ (0,\infty)
\times(0, \infty)^n $ and satisfies on this domain the
fully nonlinear partial differential inequality (PDI)
%
\begin{eqnarray}
\label{2.2}
\frac{\partial U}{\partial\tau} (\tau, \mathbf{z})
\ge\sum_{i=1}^n
\sum_{j=1}^n z_i z_j a_{ij} \biggl( \frac{1}{ 2 } D^2_{ij} U(\tau,
\mathbf{z})+ \frac{D_{i} U(\tau, \mathbf{z})}{z_1
+ \cdots+z_n} \biggr) \nonumber\\[-8pt]\\[-8pt]
\eqntext{\forall a \in\mathcal{A} (\mathbf{ z}) .}
\end{eqnarray}
We shall denote by $ \mathcal{ U} $ the collection of all
such continuous functions $ U \dvtx[0,\infty) \times(0, \infty)^n
\rightarrow(0, \infty) $ which are of class $ \mathcal{C}^{1,2} $ on
$ (0,\infty) \times(0, \infty)^n $ and satisfy (\ref{2.1})
and (\ref{2.2}). This collection $ \mathcal{ U} $ is nonempty,
since we can take $ U(\cdot, \cdot)
\equiv1 $; however, $ \mathcal{ U} $ need not contain only
one element.

Let us fix an initial configuration $ \mathbf{x} \in(0,
\infty)^n $ and consider any admissible system $ \mathcal{M} \in
\mathfrak{M} (\mathbf{x}) $. Applying It\^{o}'s rule to the process
%
\begin{equation}
\label{2.4}
\Xi(t) := X(t) L(t) U \bigl(T-t, \mathfrak{X} (t) \bigr) ,\qquad
0 \le t \le T
\end{equation}
in conjunction with (\ref{1.c}) and (\ref{1.1}), we obtain its $
\mathbb{P}^\mathcal{M}$-semimartingale decomposition as
\begin{eqnarray}
\label{2.5}
\frac{\mathrm{d} \Xi(t)}{ X(t) L(t) } &=& \Delta(t, \mathfrak{X})
\,\mathrm{d} t \nonumber\\[-8pt]\\[-8pt]
&&{}+
\sum_{\nu=1}^n \bigl[ R_\nu(t, \mathfrak{X}) - U \bigl(T-t,
\mathfrak{X}(t) \bigr) \widetilde{\vartheta}_\nu(t, \mathfrak{X}) \bigr]
\,\mathrm{d}W_\nu(t) .\nonumber
\end{eqnarray}
Here we have used the notation of (\ref{1.d}), and have set
%
\begin{eqnarray}
\label{2.6}
\qquad R_\nu(t, \mathfrak{X})&:=& \sum_{i=1}^n X_i (t) D_i
U \bigl(T-t,\mathfrak{X}(t) \bigr) \sigma_{i
\nu}(t, \mathfrak{X}) ,
\\
\label{2.7}
\Delta(t, \mathfrak{X}) &:=&
\frac{1}{ 2 }
\sum_{i=1}^n \sum_{j=1}^n X_i (t) X_j (t) \alpha_{ij}(t, \mathfrak{X})
D^2_{ij} U \bigl(T-t, \mathfrak{X}(t) \bigr)\nonumber
\\
&&{}+\sum_{i=1}^n \Biggl( \sum_{j=1}^n \frac{X_j (t)
\alpha_{ij}(t, \mathfrak{X})}{ X_1 (t)
+ \cdots+ X_n (t) } \Biggr) X_i (t)
D_{i} U \bigl(T-t, \mathfrak{X}(t) \bigr)\\
&&{} -\frac{\partial U}{\partial
\tau} \bigl(T-t, \mathfrak{X}(t) \bigr) .\nonumber
\end{eqnarray}

From the inequality of (\ref{2.2}), coupled with the fact that
$ \alpha(t, \mathfrak{X}) \in\mathcal{A} (\mathfrak{X}(t)) $ holds
for all
$ 0 \le t < \infty$, this last expression is clearly
not positive. As a~result, the positive process $ \Xi
(\cdot) $ of (\ref{2.4}) is a $
\mathbb{P}^\mathcal{M}$-supermartingale, namely,
\begin{eqnarray}
\label{2.7.7}
L (t) X(t) U \bigl(T-t, \mathfrak{X} (t) \bigr) &=& \Xi(t) \ge
\mathbb{E}^{ \mathbb{P}^\mathcal{M}}
[ \Xi(T) | \mathcal{F}(t) ] \nonumber\\[-8pt]\\[-8pt]
&=& \mathbb{E}^{ \mathbb{P}^\mathcal{M}} [
L(T)X(T) | \mathcal{F}(t) ]\nonumber
\end{eqnarray}
holds $ \mathbb{P}^\mathcal{M}$-a.s., $ \forall\mathcal{M} \in
\mathfrak{M}(\mathbf{x}) $ and $ 0 \le t \le T $; in particular,
\begin{eqnarray*}
(x_1+ \cdots+x_n) U(T, \mathbf{x}) &=& \Xi(0) \ge
\mathbb{E}^{ \mathbb{P}^\mathcal{M}}
[ \Xi(T) ] \\
&=& \mathbb{E}^{ \mathbb{P}^\mathcal{M}} [ L(T)
X(T) ] \qquad  \forall\mathcal{M} \in\mathfrak{M}
(\mathbf{x}) .
\end{eqnarray*}
With the notation of (\ref{1.8}), we obtain in this manner the following
analog of the inequality in Proposition \ref{Proposition_1}:
%
\begin{equation}
\label{2.7.b}
U (T,\mathbf{x}) \ge\Phi(T,\mathbf{x}) .
\end{equation}
%
Digging in this same spot, just a bit deeper, leads to our next result;
this is very much in the spirit of Theorem 5 in \citet{FleVer89} and of Section II.2 in Lions (\citeyear{Lio84}).

\begin{prop}
\label{Proposition_2}
For \textup{every} horizon $ T \in(0, \infty) $, initial configuration $
\mathbf{ x} \in(0, \infty)^n $ and
function $ U \dvtx[0,\infty) \times(0, \infty)^n
\rightarrow(0, \infty) $ in the collection $ \mathcal{U} $,
we have the inequality
%
\begin{equation}
\label{2.8}
U (T,\mathbf{x}) \ge\mathfrak{u} (T,\mathbf{x}) \ge\Phi
(T,\mathbf{x}) .
\end{equation}
Furthermore, the Markovian investment rule $ \pi^U \in\mathfrak{P}_*
$ generated by this function $ U $ through
\begin{eqnarray}
\label{2.11.a}
\pi^U_i (t, \mathbf{z}) &:=& z_i D_i \log
U (T-t,\mathbf{z} )\nonumber\\[-8pt]\\[-8pt]
&&{}+ \frac{z_i }{ z_1 + \cdots+
z_n } ,\qquad  (t,\mathbf{z}) \in[0,T] \times(0,
\infty)^n\nonumber
\end{eqnarray}
for each $ i = 1, \ldots, n $, satisfies for every admissible system $
\mathcal{M} \in\mathfrak{M}(\mathbf{x}) $ the inequality
%
\begin{equation}
\label{2.11.b}
Z^{ U(T, \mathbf{x}) X(0), \pi^U} ( T) \ge X(T) ,\qquad
\mathbb{P}^\mathcal{M}\mbox{-a.s.}
\end{equation}
\end{prop}

\begin{pf}
For a fixed initial configuration $
\mathbf{x} \in(0, \infty)^n $, an arbitrary admissible model
$ \mathcal{M} \in\mathfrak{M} (\mathbf{x}) $ and any function $ U
\in
\mathcal{U} $, let us recall the notation of~(\ref{2.4}) and re-cast
the dynamics of~(\ref{2.5}) as
%
\begin{equation}
\label{2.9}
\mathrm{d} \Xi(t) = \Xi(t) \Biggl( \sum_{\nu=1}^n \Psi_\nu
(t, \mathfrak{X})\, \mathrm{d} W_\nu(t) - \mathrm{d} C(t ) \Biggr) .
\end{equation}
Here by virtue of (\ref{2.5}), (\ref{2.6}) and (\ref{1.d}) we have written
\begin{eqnarray}
\label{2.10}
\Psi_\nu(t, \mathfrak{X}) &:=& \sum_{i=1}^n \sigma_{i \nu}(t,
\mathfrak
{X}) \biggl( X_i (t)
D_i \log U \bigl(T-t,\mathfrak{X}(t) \bigr) + \frac{X_i (t)}{ X
(t) } \biggr)\nonumber\\[-8pt]\\[-8pt]
&&{}- \vartheta_\nu(t, \mathfrak{X})\nonumber
\end{eqnarray}
for $ \nu=1, \ldots,n $ and have introduced in the notation of (\ref
{2.7}) the continuous, increasing process
%
\begin{equation}
\label{2.11}
C(t ) := \int_0^t \frac{ ( -\Delta(s, \mathfrak{X}) ) }{
U (T-s, \mathfrak{X}(s) ) } \,\mathrm{d} s ,\qquad  0 \le t \le T .
\end{equation}

The expression of (\ref{2.10}) suggests considering the Markovian
investment rule $ \pi^U \in\mathfrak{P}_* $ as in
(\ref{2.11.a}); then we cast the expression of (\ref{2.10}) as
\[
\Psi_\nu(t, \mathfrak{X}) = \sum_{i=1}^n \sigma_{i \nu}(t,
\mathfrak
{X}) \pi^U (t,\mathfrak{X}(t) ) - \vartheta_\nu(t, \mathfrak{X}) .
\]
On the strength of (\ref{1.6}), the value process generated by this
investment rule $ \pi^U $ starting with initial wealth $ \xi:= U(T,
\mathbf{x}) X(0) \equiv\Xi(0) $, satisfies the equation
\[
\mathrm{d} ( L( t) Z^{ \xi, \pi^U} ( t) ) = ( L( t)
Z^{ \xi, \pi^U} ( t) ) \sum_{\nu=1}^n \Psi_\nu(t, \mathfrak{X})
\,\mathrm{d}W_\nu(t) .
\]
Juxtaposing this to (\ref{2.9}), and using the positivity of $
\Xi(\cdot) $ along with the nonnegativity and nondecrease of $
C(\cdot)$, we obtain the $ \mathbb{P}^\mathcal{M}$-a.s.
comparison $ L( \cdot) Z^{ \xi, \pi^U} ( \cdot) \ge\Xi(\cdot) $, thus
%
\begin{equation}
\label{2.11.bbb}
Z^{ \xi, \pi^U} ( t) \ge X(t)
U \bigl(T-t, \mathfrak{X} (t) \bigr) ,\qquad  0 \le t \le T .
\end{equation}
With $ t = T $ this leads to (\ref{2.11.b}), in conjunction
with (\ref{2.1}). We conclude from~(\ref{2.11.b}) that the number $ U(T,
\mathbf{x})>0 $ belongs to the set on the right-hand side of
(\ref{1.7}), and the first comparison in (\ref{2.8}) follows; the
second is just a~restatement of (\ref{1.8}).
\end{pf}

\begin{cor*}
Suppose that the function $ \Phi(
\cdot, \cdot)$ of (\ref{1.8}) belongs to the collection $ \mathcal{ U}
$. Then $ \Phi( \cdot, \cdot) $ is the smallest element of $
\mathcal
{ U} $; the infimum in (\ref{1.7}) is attained; we can take $ U\equiv
\Phi$ in (\ref{2.11.b}), (\ref{2.11.a}); and the inequality in (\ref
{2.8}) holds as equality, that is, $ \Phi( \cdot, \cdot)$ coincides
with the arbitrage function
\[
\mathfrak{u} (T,\mathbf{x}) = \Phi(T,\mathbf{x}) \qquad  \forall
(T,\mathbf
{x}) \in(0,\infty)
\times(0,\infty)^n .
\]
\end{cor*}

\textit{Interpretation:} Imagine that the small investor
is a manager who invests for a pension fund and tries to track or
exceed the performance of an index (the market portfolio) over a finite
time-horizon. He has to do this in the face of uncertainty
about the characteristics of the market, including its covariance
and price-of-risk structure, so he acts with extreme prudence and
tries to protect his clients against the most adverse market
configurations imaginable [the range of such configurations is
captured by the constraints~(\ref{E.3}), (\ref{1.b})]. If such adverse
circumstances do not materialize his strategy generates a
surplus, captured here by the increasing process $ C(\cdot) $
of (\ref{2.11}) with $ U \equiv\Phi\equiv\mathfrak{u} $, which can
then be returned to the (participants in the) fund. We are borrowing
and adapting this interpretation from \citet{LYO95}.

Similarly, the Markovian investment rule $ \pi^U \in\mathfrak{P}_* $
generated by the function $ U \equiv\Phi\equiv\mathfrak{u} $ in
(\ref
{2.11.a}), (\ref{FG}) implements the best possible outperformance of
the market portfolio, as in (\ref{2.11.b}).


\section{A verification result}
\label{sec6}

For the purposes of this section we shall impose the following \textit
{growth} condition on the family $ \mathbb{A}= \{ \mathcal{A}
(\mathbf
{y}) \}_{\mathbf{y} \in\mathfrak{S}_n} $ of subsets of~$ \mathbb{S}^n
$ in (\ref{E.3.a}), (\ref{1.b}): there exists a constant
$ C \in( 0, \infty) $, such that for all $ \mathbf{y} \in\mathfrak
{S}_n $ we have
%
\begin{equation}
\label{3.2.GG}
\sup_{a \in\mathcal{A} (\mathbf{y})} \biggl( \max_{1 \le i,j \le
n} { y_i y_j |a_{ij}| \over(y_1 + \cdots
+ y_n ) } \biggr) \le C ( 1 + \|\mathbf{ y}\| ) .
\end{equation}
We shall also need the following \textit{strong
ellipticity} condition, which mandates that for every nonempty,
compact subset $ \mathbf{K} $ of $ (0, \infty)^n $, there
exists a real constant $ \lambda= \lambda_\mathbf{ K} > 0 $ such that
%
\begin{equation}
\label{3.2}
\inf_{\mathbf{ z} \in\mathbf{ K}} \Biggl( \inf_{a \in\mathcal{A}
(\mathbf
{z})} \Biggl( \sum_{i=1}^n \sum_{j=1}^n \xi_i \xi_j
a_{ij} \Biggr) \Biggr) \ge\lambda_\mathbf{ K} \|\xi\|^2 \qquad  \forall
\xi\in\mathbb{R}^n .
\end{equation}

\renewcommand{\theass}{\Alph{ass}}
\begin{ass}\label{assA}
There exist a continuous
function $ \mathbf{a}\dvtx(0, \infty)^n \rightarrow\mathbb{S}^n $,
a $ \mathcal{C}^2$-function $ H \dvtx(0, \infty)^n \rightarrow
\mathbb{R} $, and a continuous square root $ \mathbf{ s}
(\cdot) $ of $ \mathbf{a} (\cdot) $, namely $ \mathbf{a}
(\cdot) = \mathbf{ s} (\cdot) \mathbf{ s}' (\cdot) $ such
that, with the vector-valued function $ \theb(\cdot) = ( \theb_1
(\cdot), \ldots, \theb_n (\cdot) )' $ defined by
%
\begin{equation}
\label{D.2}
\theb_\nu( \mathbf{ z}) := \sum_{j=1}^n z_j \mathbf{s}_{j \nu
}(\mathbf
{ z}) D_j H(\mathbf{ z}) ,\qquad  \nu= 1, \ldots, n ,
\end{equation}
condition (\ref{D.22}) is satisfied, whereas the system of stochastic
differential equations
%
\begin{eqnarray}
\label{D.1}
\mathrm{d} X_i(t) = X_i(t)
\sum_{\nu=1}^n \mathbf{s}_{i \nu}( \mathfrak{ X} (t) ) [ \mathrm{d}
W_\nu(t) +\theb_\nu( \mathfrak{ X} (t)) \,\mathrm{d} t ] ,\nonumber\\[-8pt]\\[-8pt]
\eqntext{X_i(0)=x_i>0, i=1,\ldots,n}
\end{eqnarray}
has a solution in which the state process $ \mathfrak{ X} (\cdot) $
takes values in $ (0, \infty)^n $.
\end{ass}

A bit more precisely, this assumption posits the existence of
a Markovian admissible system $ \mathcal{M}_o \in\mathfrak{ M}_*
(\mathbf{
x}) $ consisting of a filtered probability space $
(\Omega,\mathcal{F}, \mathbb{P})$, $ \mathbb{F} = \{
\mathcal{F}(t)\}_{0 \le t <\infty} $ and of two continuous,
adapted process $ \mathfrak{X}(\cdot) $ and~$ W(\cdot) $ on
it, such that under the probability measure~$ \mathbb{P} \equiv
\mathbb{P}^{\mathcal{M}_o} $ the process~$ W(\cdot) $ is
$ n$-dimensional Brownian Motion, the process
$ \mathfrak{X}(\cdot) $ takes values in $ (0, \infty)^n $
a.s. and (\ref{1.1}) holds with $ \vartheta_\nu(t, \mathfrak{
X}) = \theb_\nu( \mathfrak{ X} (t)) $ as in (\ref{D.2}), and with $
\sigma_{i \nu}
(t,\mathfrak{ X})= \mathbf{s}_{i \nu}( \mathfrak{ X} (t) ) $, $0\le
t <
\infty$ ($1 \le i, \nu\le n$). The system of equations (\ref{D.1})
can be cast equivalently as
\begin{eqnarray}
\label{D.5}
\mathrm{d} X_i(t) &=& X_i(t) \Biggl[
\sum_{\nu=1}^n \mathbf{s}_{i \nu}( \mathfrak{ X} (t) ) \,\mathrm{d}
W_\nu(t)\nonumber\\[-8pt]\\[-8pt]
&&\hphantom{X_i(t) \Biggl[}
{} + \Biggl( \sum_{j=1}^n
\mathbf{a}_{i j}( \mathfrak{ X} (t) ) X_j (t) D_j H (\mathfrak{ X} (t))
\Biggr)\, \mathrm{d} t \Biggr].\nonumber
\end{eqnarray}

\begin{ass}\label{assB}
In the notation of the previous paragraph and
under the condition
%
\begin{equation}
\label{D.4}
\sum_{i=1}^n \sum_{\nu=1}^n z_i |\mathrm{s}_{i \nu}(\mathbf{ z})
| |
\theb_\nu(\mathbf{ z}) | \le C (1+ \| \mathbf{ z}\| ) \qquad  \forall
\mathbf{ z} \in(0, \infty)^n ,
\end{equation}
we define on $ (0, \infty)^n $ the continuous functions
$
g(\mathbf{ z}) :=e^{ - H(\mathbf{ z})} \sum_{i=1}^n
z_i $ and $ k(\mathbf{z}) := (1/2) \sum_{i=1}^n \sum_{j=1}^n
\mathbf{a}_{ij} (\mathbf{z}) [ D^2_{i j}
H(\mathbf{ z}) + D_i H(\mathbf{ z}) D_j H(\mathbf{ z}) ] $ and assume
that the function
\[
G (\tau, \mathbf{ x}) := \mathbb{E}^{\mathbb{P}^{\mathcal{M}_o}} \Biggl[
g (
\mathfrak{ X} (\tau) ) \exp\Biggl\{ \int_0^\tau k ( \mathfrak{ X} (t) )
\,\mathrm{d} t \Biggr\} \Biggr] ,\qquad  (\tau, \mathbf{x}) \in[0, \infty)^n
\]
is continuous on $ [0, \infty) \times(0, \infty) $ and of class $
\mathcal{C}^{1,2} $ on $ (0, \infty) \times(0, \infty) .$
\end{ass}

Sufficient conditions for Assumptions \ref{assA}, \ref{assB} to hold are provided in
\citet{FerKar10}, Sections 8 and 9. It is also shown there,
that we have the $ \mathbb{P}^{\mathcal{M}_o}$-martingale property
%
\begin{equation}
\label{D.8}
\quad \mathbb{E}^{\mathbb{P}^{\mathcal{M}_o}} [ X(T) L(T) | \mathcal{F} (t)
] = X(t) L(t) \cdot\Gamma\bigl(T-t, X (t) \bigr) ,\qquad  0 \le t \le T
\end{equation}
for the function
%
\begin{equation}
\label{D.9} \Gamma(\tau, \mathbf{ z}) := G (\tau, \mathbf{z})
/ g ( \mathbf{ z}) ,\qquad  (\tau, \mathbf{ z}) \in[0, \infty)
\times(0, \infty)^n .
\end{equation}
This function is of class $ \mathcal{C}^{1,2} $ on $ (0, \infty)
\times
(0, \infty) $ and satisfies the initial condition $ \Gamma(0, \cdot)
\equiv1 $ on $ (0, \infty)^n $ as well as the \textit{linear}
second-order parabolic equation
%
\begin{eqnarray}
\label{D.10}
\frac{\partial\Gamma}{\partial\tau} (\tau, \mathbf{z}) = \sum_{i=1}^n
\sum_{j=1}^n z_i z_j \mathbf{a}_{ij} (\mathbf{ z}) \biggl( \frac{1}{ 2 }
D^2_{ij} \Gamma(\tau,\mathbf{z})+ \frac{D_{i} \Gamma(\tau,
\mathbf
{z})}{z_1+ \cdots+z_n} \biggr) ,\nonumber\\[-8pt]\\[-8pt]
\eqntext{(\tau, \mathbf{ z}) \in(0, \infty)
\times(0, \infty)^n .}
\end{eqnarray}

\begin{prop}[(Verification argument)]
\label{Proposition_3}
Under the Assumptions \ref{assA}, \ref{assB} and the
conditions (\ref{3.2.GG}), (\ref{3.2}), suppose that the functions $
\mathbf{a} (\cdot) $ and $ \Gamma(\tau, \cdot) $ satisfy the inequality
\begin{eqnarray}
\label{D.11}
\qquad &&\sum_{i=1}^n
\sum_{j=1}^n z_i z_j \mathbf{a}_{ij} (\mathbf{ z}) \biggl( \frac{1}{ 2 }
D^2_{ij} \Gamma(\tau,\mathbf{z})+ \frac{D_{i} \Gamma(\tau,
\mathbf
{z})}{z_1+ \cdots+z_n} \biggr)
\nonumber\\[-8pt]\\[-8pt]
\qquad &&\qquad \ge\sum_{i=1}^n
\sum_{j=1}^n z_i z_j a_{ij} \biggl( \frac{1}{ 2 } D^2_{ij} \Gamma(\tau
,\mathbf{z})+ \frac{D_{i} \Gamma(\tau, \mathbf{z})}{z_1+ \cdots+z_n}
\biggr) \qquad  \forall a \in\mathcal{A} (\mathbf{ z})\nonumber
\end{eqnarray}
for every $ (\tau, \mathbf{ z}) \in(0, \infty) \times(0, \infty
)^n $.
Then, in the notation of (\ref{1.7})--(\ref{1.9}), we have:
\[
\mathfrak{u} (T,\mathbf{x}) = \Phi(T,\mathbf{x}) =
\Gamma(T,\mathbf{x}) = \mathfrak{u}_{\mathcal{M}_o} (T,\mathbf{x}) ,\qquad
\forall(T,\mathbf{x}) \in
(0,\infty) \times(0,\infty)^n
\]
for the Markovian admissible system $ \mathcal{M} \equiv\mathcal{M}_o
\in\mathfrak{M} (\mathbf{ x}) $ posited in Assumption~\ref{assA}; the
conclusions of Proposition \ref{Proposition_2} and its corollary for $
U \equiv\Phi$; as well as the $ \mathbb{P}^\mathcal{M}$-a.s. comparison
\[
L (t) X(t) \cdot\mathfrak{u} \bigl(T-t, \mathfrak{X} (t) \bigr) \ge
\mathbb{E}^{ \mathbb{P}^\mathcal{M}} [ L(T)X(T) | \mathcal{F}(t) ],\qquad  0\le t \le T ,
\]
which holds for every $ \mathcal{M} \in\mathfrak{ M} (\mathbf{ x}) $
and as equality for $ \mathcal{M} \equiv\mathcal{M}_o \in\mathfrak{M}
(\mathbf{ x}) $.
\end{prop}

\begin{pf}
Under condition (\ref{D.11}) the
function $ \Gamma(\cdot, \cdot) $ belongs to the collection~$ \mathcal{ U} $, as (\ref{2.2}) is satisfied with $ U \equiv\Gamma$
on the strength of (\ref{D.10}) and (\ref{D.11}); thus, we deduce $
\Gamma(T, \mathbf{x}) \ge\Phi(T, \mathbf{x}) $ from (\ref{2.7.b}).
On the other hand, equality~(\ref{D.8}) with $ t=0 $, and the
definition of $ \Phi(T, \mathbf{x}) $ in (\ref{1.8}), give
\[
\Gamma(T,\mathbf{x}) =
\frac{ \mathbb{E}^{ \mathbb{P}^{\mathcal{M}_o}}
[ L(T) X(T) ] }{ x_1 + \cdots+ x_n } = \mathfrak{u}_{\mathcal{M}_o}
(T,\mathbf{x}) \le\Phi(T,\mathbf{x}) ,
\]
so the equality $ \Gamma(T, \mathbf{x}) = \Phi(T,
\mathbf{x}) $ follows. In other words, we identify $ \mathcal{M}_o $
as a
Markovian admissible system that satisfies (\ref{1.88}) and attains the
supremum in (\ref{1.8}). The remaining claims come from Proposition
\ref{Proposition_2} and its Corollary [in particular, from reading (\ref
{2.8}) with $ U \equiv
\Gamma$] and from (\ref{2.7.7}), (\ref{D.8}).
\end{pf}

\begin{remark}
\label{Remark_2}
Proposition \ref{Proposition_3} holds under conditions weaker
than those imposed in Assumptions \ref{assA} and \ref{assB} above, at the ``expense'' of a
certain localization. More precisely, one posits the existence of
locally bounded and locally Lipschitz functions $ \mathbf{s}_{i \nu}
(\cdot) $ and $ \theb_\nu(\cdot) $ ($1 \le i, \nu\le n$) for which
(\ref{D.22}), (\ref{3.2.GG}), (\ref{3.2}) and (\ref{D.4}) are satisfied
with $ \mathbf{a} (\cdot) = \mathbf{ s} (\cdot) \mathbf{ s}'
(\cdot) $,
and for which there exists a Markovian admissible system $ \mathcal
{M}_o \in\mathfrak{ M}_* (\mathbf{x}) $ whose state process $
\mathfrak
{ X} (\cdot) $ in (\ref{D.1}) is, under $ \mathbb{P}^{ \mathcal{M}_o}$,
a strong Markov process with values in $ (0, \infty)^n$ a.s. Using
results from the theory of stochastic flows [\citet{Kun90}, \citet{Pro04}]
and from parabolic partial differential equations [\citet{JanTys06}, Ekstr\"{o}m and Tysk (\citeyear{EksTys09})], Theorem~2 in Ruf (\citeyear{RUF})
shows that the function $ \Gamma(\cdot,\cdot) $ is then of class $
\mathcal{C}^{1,2} $ locally on $ (0,\infty) \times(0,\infty)^n $, and
solves there equation (\ref{D.10}).
\end{remark}


\section{Maximizing the probability of containment}
\label{sec7}

We have now gone as far as we could without having to specify the
nature of our filtered measurable space $ (\Omega, \mathcal{F} )$, $
\mathbb{F} = \{ \mathcal{F}(t)\}_{0 \le t <\infty} $. To proceed
further, we shall need to choose this space carefully.

We shall take as our sample space the set $ \Omega$ of
right-continuous paths $ \omega\dvtx[0,\infty) \rightarrow\mathfrak{S}_n
\cup\{ \Delta\} $. Here $ \Delta$ is an additional ``absorbing
point''; paths stay at $ \Delta$ once they get there, that is, after
$ \mathcal{T} (\omega) = \inf\{ t \ge0 | \omega(t) =
\Delta\} $, and are continuous on $ (0, \mathcal{T} (\omega))$; we
are employing here, and throughout this work, the usual convention $
\inf\varnothing= \infty$. We also select $ \mathcal{K} (\Delta) = \{
(\mathbf{ 0},
\mathrm{O}_{n \times n}) \} $ and $ \mathcal{A} (\Delta) =
\{ \mathrm{O}_{n \times n} \} $, where
$ \mathrm{O}_{n \times n} $ is the zero matrix. With
$ \mathcal{F}^{ \flat} (t) :=\sigma(\omega(s), 0\le s \le t) $,
the filtration $ \mathbb{F}^{ \flat} = \{ \mathcal{F}^{ \flat} (t)\}
_{0 \le t < \infty}
$ is a \textit{standard system} in the terminology of \citet{Par67}. This means that each $ ( \Omega, \mathcal{F}^{ \flat} (t) ) $
is isomorphic to the Borel $\sigma$-algebra on some Polish space, and
that for any decreasing sequence $ \{ A_j\}_{j \in\mathbb{N}} $ where
each $A_j$ is an atom of the corresponding $ \mathcal{F}^{ \flat} (t_j)
$, for some increasing sequence $ \{ t_j\}_{j \in\mathbb{N}} \subset
[0, \infty) $, we have $ \bigcap_{j \in\mathbb{N}} A_j \neq\varnothing$
[see the Appendix in F\"{o}llmer (\citeyear{Fol72}), as well as \citet{Mey72} and F\"{o}llmer (\citeyear{Fol73})].

With all this in place we take $ ( \Omega, \mathcal{F} ), \mathbb{F} =
\{
\mathcal{F} (t)\}_{0 \le t < \infty} $ as our filtered measurable
space, where
\[
\mathcal{F} (t) := \bigcap_{\varepsilon
>0} \mathcal{F}^{ \flat} (t+\varepsilon)\quad  \mbox{and}\quad
\mathcal{F} := \sigma\biggl( \bigcup_{0 \le t < \infty}
\mathcal{F} (t ) \biggr) .
\]

An admissible system $ \mathcal{M} \in\mathfrak{ M} (\mathbf{
x}) $, $\mathbf{ x} \in(0, \infty)^n $ defined as in Section \ref{sec2}
consists of this filtered measurable space $ ( \Omega, \mathcal{F} ),
\mathbb{F} = \{\mathcal{F} (t)\}_{0 \le t < \infty} $, of a probability
measure $ \mathbb{P}^\mathcal{M} $ on it, of an $ n$-dimensional
Brownian motion $ W(\cdot) $ on the resulting probability space and of
the coordinate mapping process $ \mathfrak{ X}(t,\omega) = \omega(t) ,
0 \le t < \infty$ which is assumed to satisfy
(\ref{1.1}), (\ref{1.2}) and to take values in $ (0, \infty)^n $, $
\mathbb{P}$-a.s. We shall take
\[
\mathcal{ T} :=
\inf\{ t \ge0 | \Lambda(t )=0\} = \inf\{ t \ge0 | L (t ) X(t
)=\infty\}
\]
in the notation of (\ref{Lambda}), (\ref{expo}) and (\ref{1.a}), and
note $ \mathbb{P}^\mathcal{M} ( \mathcal{T} < \infty) =0 $.

\subsection{The F\"{o}llmer exit measure}
\label{subsec7.1}

With this setup, there exists a probability measure $
\mathbb{Q} $ on $ ( \Omega, \mathcal{F} ), $ such that
%
\begin{equation}
\label{bigtriangleleft}
\mathrm{d} \mathbb{P}^\mathcal{M} = \Lambda(T) \,\mathrm{d} \mathbb{Q}
\qquad\mbox{holds on each } \mathcal{F} (T) ,\qquad  T \in(0, \infty) ;
\end{equation}
we express this property (\ref{bigtriangleleft}) by writing $ \mathbb
{P}^\mathcal{M} \ll\mathbb{Q} $. Under the measure~$ \mathbb{Q} $,
the process $ \widetilde{W} (\cdot) $ of (\ref{1.1.b}) is Brownian
motion; whereas the processes $ \mu_1 (\cdot), \ldots,\allowbreak\mu_n
(\cdot) $ and $ \Lambda(\cdot) $ of (\ref{MU}),
(\ref{Lambda}) in Section \ref{subsec3.1} are nonnegative $
\mathbb{Q}$-martingales.

The ``absorbing state'' $ \Delta$ acts here as a proxy for $
\mathbb{P}^\mathcal{M}$-null sets to which the new measure $ \mathbb{Q} $ may assign
positive mass; the possible existence of such sets makes it important
that the filtration $ \mathbb{F} $ be ``pure,'' that is, \textit{not}
completed by $ \mathbb{P}^\mathcal{M}$-null sets. This probability measure $
\mathbb
{Q} $ satisfies
\begin{eqnarray}
\label{Foell}
\frac{ \mathbb{E}^{\mathbb{P}^\mathcal{M}}[ L(T) X(T) ] }{ x_1 +
\cdots
+ x_n } &=& \mathbb{E}^{\mathbb{P}^\mathcal{M}} \bigl[ ( 1 / \Lambda(T) )
\cdot\mathbf{ 1}_{\{ \mathcal{T} >T\}} \bigr] \nonumber\\[-8pt]\\[-8pt]
&=& \mathbb{Q} ( \mathcal
{T} >T) \qquad  \forall T \in[0, \infty)\nonumber
\end{eqnarray}
and
%
\begin{equation}
\label{Ttau.a}
\mathcal{ T} = \inf\Biggl\{ t \ge0 \Big| \int_0^t \| \widetilde{\vartheta} (s,
\mathfrak{ X}) \|^2 \,\mathrm{d} s =\infty\Biggr\} ,\qquad  \mathbb{Q}\mbox{-a.s.}
\end{equation}

We also have $\mathbb{Q}$-a.e. on $ \{ T<\mathcal{T}
<\infty\} $
\[
L (\mathcal{T}+u) X(\mathcal{T}+u) = \infty\qquad  \forall u \ge0
\]
and
\[\int_0^T \| \widetilde{\vartheta} (t, \mathfrak{X}) \|^2
\,\mathrm{d} t < \int_0^\mathcal{T} \| \widetilde{\vartheta} (t, \mathfrak
{X}) \|^2 \,\mathrm{d} t = \infty.
\]
Whereas, $\mathbb{Q} $-a.e. on $ \{ \mathcal{T} =\infty\}
$, we have
\[
L (T) X(T) < \infty,\qquad  \int_0^T \| \widetilde{\vartheta} (t,
\mathfrak
{X}) \|^2 \,\mathrm{d} t < \infty;\qquad  \forall T \in[0, \infty) .
\]

We deduce from (\ref{Foell}) that the arbitrage function of (\ref{1.9})
for the model $ \mathcal{M} \in\mathfrak{ M} (\mathbf{ x}) $ is given
by the \textit{probability of} ``\textit{containment}'' \textit{under the measure}~$
\mathbb
{Q} $, namely, the $ \mathbb{Q}$-probability that the process $
\mathfrak{ X}(\cdot) $, started at $ \mathbf{ x} \in(0, \infty)^n $,
stays in $ (0, \infty)^n $ throughout the time-horizon $ [0,T] $.

At this point we shall impose the following requirements on $ \mathbb
{K}= \{ \mathcal{K} (\mathbf{y}) \}_{\mathbf{y} \in\mathfrak{S}_n} $,
the family of compact, convex
subsets of $ \mathbb{R}^n \times\mathbb{S}^n $ in (\ref{1.b}),
(\ref
{E.3}): there exists a constant $ 0< C < \infty
$, such that for all $ \mathbf{y} \in\mathfrak{S}_n $ we have the
strengthening
%
\begin{equation}
\label{3.2.G}
\sup_{a \in\mathcal{A} (\mathbf{y})} \Biggl( \sum_{i=1}^n \sum_{i=1}^n
y_i y_j a_{ij} \Biggr) \le C (y_1 + \cdots+ y_n )^2
\end{equation}
of the growth condition in (\ref{3.2.GG}), as well as the ``shear'' condition
%
\begin{equation}
\label{ThetaTrace}
\sup_{(\theta, a) \in\mathcal{K} (\mathbf{ y})} \biggl[
\biggl( { \| \theta\|^2 \over1 + \operatorname{Tr} (a) } \biggr) +
\biggl( { \operatorname{Tr} (a) \over1 + \| \theta\|^2
} \biggr) \biggr] \le C .
\end{equation}
Then the following identity holds $\mathbb{Q} $-a.s.:
%
\begin{equation}
\label{Ttaui}
\mathcal{ T} = \min_{1 \le i \le n} \mathcal{ T}_i \qquad  \mbox{where }
\mathcal{ T}_i := \inf\{ t \ge0 | X_i (t)=0\} .
\end{equation}
For justification of the claims made in this subsection, we refer to
Section~7 in Fernholz and Karatzas (\citeyear{FerKar10}), as well as \citet{DelSch95N2}, \citet{PALPRO10} and \citet{RUF}---in
addition, of course, to the seminal work by F\"{o}llmer (\citeyear{Fol72},\vadjust{\goodbreak} \citeyear{Fol73}).

The special
structure of the filtered measurable space $ ( \Omega, \mathcal{F} ),
\mathbb{F} = \break\{
\mathcal{F} (t)\}_{0 \le t < \infty} $ that we selected in this section
is indispensable for this construction and for the representation (\ref
{Foell}); whereas the inequality $ \| \theta\|^2 \le C ( 1 + \mathrm
{Tr} (a)), \forall(\theta,
a) \in\mathcal{K} (\mathbf{ y}) , \mathbf{ y} \in\mathfrak{S}_n $
from condition (\ref{ThetaTrace}) is important for establishing the
representation of (\ref{Ttaui}).

\subsection{Auxiliary admissible systems}
\label{subsec7.2}

Let us fix then an initial configuration $ \mathbf{ x} = (x_1, \ldots,
x_n )' \in(0, \infty)^n $ and denote by $ \mathfrak{N} (\mathbf{
x}) $
the collection of stochastic systems $ \mathcal{N} $ that consist of
the filtered measurable space $ (\Omega, \mathcal{F} ), \mathbb{F} =
\{
\mathcal{F} (t) \}_{0 \le t < \infty} $, of a probability measure $
\mathbb{Q} \equiv\mathbb{Q}^{ \mathcal{N}} $,
of an $ \mathbb{R}^n$-valued Brownian motion $ \mathbb{\widetilde{W}}
(\cdot) $ under $ \mathbb{Q} $, and of the co\"{o}rdinate mapping
process $ \mathfrak{ X}(t, \omega) = \omega(t) , (t, \omega) \in[0,
\infty) \times\Omega$ which satisfies $ \mathbb{Q} $-a.s. the
system of the stochastic equations in (\ref{1.1.a})
%
\begin{eqnarray}
\label{aux}
\mathrm{d} X_i(t) &=& X_i(t)
\Biggl( { \sum_{j=1}^n \alpha_{ij} (t, \mathfrak{ X} ) X_j (t) \over X_1
(t) + \cdots+ X_n (t) } \,\mathrm{d} t + \sum_{\nu=1}^n \sigma_{i
\nu
}(t, \mathfrak{ X}) \,\mathrm{d} \widetilde{W}_\nu(t) \Biggr)
\nonumber\\
&=& X_i(t) \sum_{\nu=1}^n \sigma_{i \nu} (t, \mathfrak{ X}) \Biggl(
\,\mathrm{d}
\widetilde{W}_\nu(t) + \sum_{j=1}^n {\sigma_{j \nu} (t, \mathfrak{X}
) X_j (t) \over X_1 (t) + \cdots+ X_n (t) } \,\mathrm{d} t \Biggr) ,\\
\eqntext{i=1,\ldots, n .}
\end{eqnarray}
Here the elements $ \sigma_{i \nu} \dvtx[0, \infty) \times\Omega
\rightarrow\mathbb{R} $, $1 \le i, \nu\le n $ of the matrix $ \sigma
(\cdot, \cdot) = \{ \sigma_{i \nu} (\cdot, \cdot) \}_{1 \le i ,
\nu\le
n} $ are progressively measurable functionals that satisfy, in the
notation of (\ref{E.3.a}),
%
\begin{equation}
\label{const1}
\sigma(t, \omega) \sigma' (t, \omega)=: \alpha(t, \omega) \in
\mathcal{A} ( \omega(t) ) \qquad  \forall( t, \omega) \in[0, \infty)
\times\Omega.
\end{equation}

As in Section \ref{subsec2.3}, we shall denote by $ \mathfrak{N}_*
(\mathbf{ x}) $ the subcollection of $ \mathfrak{N} (\mathbf{ x}) $
that consists of \textit{Markovian} auxiliary admissible systems, namely,
those for which the equations of (\ref{aux}) are satisfied with $
\alpha(t, \mathfrak{X}) = \mathbf{a} (t, \mathfrak{X} (t)) $ and $
\sigma(t, \mathfrak{X}) = \mathbf{s} ( t, \mathfrak{X} (t)) $, $ 0
\le t < \infty$ and with measurable functions $ \mathbf{a} \dvtx[0,
\infty
) \times\mathfrak{S}_n \rightarrow\mathbb{S}^n $ and $ \mathbf{s}
\dvtx
[0, \infty) \times\mathfrak{S}_n \rightarrow\mathcal{L} (\R^n; \R^n)
$ that satisfy the condition $ \mathbf{s} (t, \mathbf{ y}) \mathbf{s}'
(t, \mathbf{ y}) =\mathbf{a} (t, \mathbf{ y}) \in\mathcal{A} (y)
$, $
\forall(t, \mathbf{ y}) \in[0, \infty) \times\mathfrak{S}_n $. We
invoke the same Markovian selection results as in Section~\ref
{subsec2.3}, to ensure that the process $ \mathfrak{X} (\cdot) $ is
strongly Markovian under any given $ \mathbb{Q}^\mathcal{N} $, $
\mathcal{N} \in\mathfrak{N}_* (x) $.

By analogy with (\ref{Ttaui}), we consider
%
\begin{equation}
\label{Ttaui1}
\mathcal{ \widehat{T}} (\omega) := \min_{1 \le i \le n} \mathcal{ T}_i
(\omega) \qquad \mbox{with } \mathcal{ T}_i (\omega) =
\inf\{ t \ge0 | \omega_i (t)=0\} .
\end{equation}
Then for every $ \omega\in\{ \widehat{\mathcal{T}} < \infty\} $ we have
%
\begin{equation}
\label{const2}
\quad\qquad \int_0^T \operatorname{Tr} ( \alpha(t, \omega) )\, \mathrm{d} t < \int
_0^{\mathcal{\widehat{T}}(\omega)} \operatorname{Tr} ( \alpha(t, \omega
) )\,
\mathrm{d} t = \infty\qquad  \forall0 \le T < \mathcal{\widehat{\mathcal
{T}}}(\omega) ;
\end{equation}
whereas $ \int_0^T \operatorname{Tr} ( \alpha(t, \omega) ) \,\mathrm{d} t <
\infty$, $ 0 \le T < \infty$ holds for every\vadjust{\goodbreak} $ \omega\in\{ \widehat
{\mathcal{T}} = \infty\} $.

\begin{remark}
\label{Remark_7}
As in Section \ref{subsec2.2}, solving the stochastic equation (\ref
{aux}) subject to condition (\ref{const1}) amounts to requiring that
the process
\[
u (t , \mathfrak{X} (t) ) - \int_0^t \biggl( { \partial u \over\partial s }
(s , \mathfrak{X} (s) )+ \widehat{\mathcal{L}} u ( \mathfrak{X} (s)
) \biggr)\,\mathrm{d} s ,\qquad  0 \le t < \infty
\]
be a local supermartingale, for every continuous $ u \dvtx(0, \infty)
\times\mathfrak{S}_n \rightarrow\mathbb{R} $ which is of class $
\mathcal{C}^{1,2} $ on $ (0, \infty) \times(0, \infty)^n $ and has
compact support; here $ \widehat{\mathcal{L}} $ is the nonlinear
second-order partial differential operator in (\ref{E.3.d}).
\end{remark}

\begin{remark}
\label{Remark_3}
The total capitalization process $ X(\cdot)= X_1(\cdot) + \cdots
+ X_n(\cdot) $ satisfies, by virtue of (\ref{aux}), the equation
\begin{eqnarray*}
d X(t) &\hspace*{3pt}=& X(t) [ \mathrm{d} \widetilde{N} (t) + \mathrm{d} \langle
\widetilde{N} \rangle(t) ] , \\
\widetilde{N}(\cdot) &:=& \sum_{\nu=1}^n
\int_0^\cdot\Biggl( \sum_{i=1}^n ( X_i (t) / X(t) ) \sigma_{i \nu} (t,
\mathfrak{ X}) \Biggr)\, \mathrm{d} \widetilde{W}_\nu(t) .
\end{eqnarray*}
Under the measure $ \mathbb{ Q} $, the process $ \widetilde{N} (\cdot)
$ is a continuous local martingale with quadratic variation
\[
\langle\widetilde{N} \rangle( t) = \sum_{i=1}^n \sum_{j=1}^n \int
_0^t X_i (s) \alpha_{i j} (s, \mathfrak{ X}) X_j (s) \bigl( X_1 (s) +
\cdots
+ X_n (s) \bigr)^{-2} \,\mathrm{d} s \le C t
\]
from (\ref{3.2.G}), so the total capitalization process
\[
X(t) = X(0) \cdot e^{ \widetilde{N}(t) + (1/2) \langle\widetilde{N}
\rangle(\cdot) } = X(0) \cdot e^{ \widetilde{B}(u) + (u/2)} |_{ u =
\langle\widetilde{N} \rangle(t)} ,\qquad  0 \le t < \infty
\]
takes values in $ (0, \infty) $, $ \mathbb{ Q}$-a.e.; here $
\widetilde
{B}(\cdot) $ is a one-dimensional $ \mathbb{ Q}$-Brownian motion. This
is in accordance with our selection of the punctured nonnegative
orthant $ \mathfrak{S}_n $ in (\ref{1.b}) as the state-space for the
process $ \mathfrak{ X} (\cdot) $ under $ \mathbb{Q} $.

Under this measure, the relative weights $ \mu_i (\cdot) = X_i (\cdot)
/ X(\cdot) $, $ i=1, \ldots, n $ are nonnegative local martingales and
supermartingales, in accordance with (\ref{MU}), and since $ \sum
_{i=1}^n \mu_i (\cdot) \equiv1 $ these processes are bounded, so they
are actually martingales. Once any one of the processes $ X_1 (\cdot) ,
\ldots, X_n (\cdot) $ [i.e., any one of the processes $ \mu_1
(\cdot)
, \ldots, \mu_n (\cdot) $] becomes zero, it stays at zero forever; of
course, not all of them can vanish at the same time.
\end{remark}

In Section \ref{subsec7.1} we started with an arbitrary
admissible system $ \mathcal{ M} \in\mathfrak{ M} (\mathbf{
x}) $ and produced an ``auxiliary'' admissible system $
\mathcal{N} \in\mathfrak{N} (\mathbf{ x}) $, for which the
property (\ref{Foell}) holds. Thus, for every $ (T, \mathbf{
x}) \in(0, \infty) \times(0, \infty)^n $ we deduce
\begin{eqnarray}
\label{3.7.f}
Q (T, \mathbf{ x}) &:=& \sup_{ \mathcal{N} \in
\mathfrak{N} (\mathbf{ x})} \mathbb{
Q}^{ \mathcal{N}} ( \mathcal{T} >T ) \nonumber\\[-8pt]\\[-8pt]
&\hspace*{3pt}\ge&
\sup_{\mathcal{ M} \in\mathfrak{ M} (\mathbf{ x})} \biggl(
\frac{ \mathbb{E}^{ \mathbb{P}^\mathcal{M}}[ L(T) X(T)
] }{ x_1 + \cdots+ x_n } \biggr)= \Phi(T, \mathbf{ x}) .\nonumber
\end{eqnarray}

\subsection{Preparatory steps}
\label{subsec7.2a}

We suppose from now onwards that, for every progressively measurable
functional $ \alpha\dvtx[0, \infty) \times\Omega\rightarrow
\mathbb
{S}^n $ which satisfies
%
\begin{equation}
\label{Alpha}
\alpha(t, \omega) \in\mathcal{A} ( \omega(t) ) \qquad \mbox{for all }
(t,\omega) \in[0, \infty) \times\Omega,
\end{equation}
we can select a progressively measurable functional $ \vartheta\dvtx[0,
\infty) \times\Omega\rightarrow\mathbb{R}^n $ with
%
\begin{equation}
\label{ThetaAlpha}
( \vartheta(t, \omega), \alpha(t, \omega) ) \in\mathcal{K} (
\omega(t) ) \qquad  \forall(t, \omega) \in[0, \infty) \times\Omega
\end{equation}
[see the ``measurable selection'' results in Chapter 7 of \citet{BerShr78}]. We introduce now the functional
%
\begin{equation}
\label{Ttau101}
\widetilde{\vartheta} (t, \omega) :=
\vartheta(t, \omega) - \sigma' (t , \omega) \omega(t) / \bigl( \omega_1
(t) + \cdots+ \omega_n (t) \bigr)
\end{equation}
as in (\ref{1.d}) and also, by analogy with (\ref{Ttau.a}), the
stopping rule
%
\begin{equation}
\label{Ttau1}
\mathcal{ T} (\omega)
:= \inf\biggl\{ t \ge0 \Big| \int_0^t \| \widetilde{\vartheta} (s, \omega)
\|
^2 \,\mathrm{d} s =\infty\biggr\}
\end{equation}
[cf. \citet{LevSko95}, where stopping rules of this type
also play very important roles in the study of arbitrage].

We recall now (\ref{const2}); on the strength of the requirement
(\ref{ThetaTrace}), this gives
%
\begin{equation}
\label{Ttau111}
\quad \int_0^T \| \vartheta(t, \omega) \|^2 \,\mathrm{d} t < \int
_0^{\mathcal
{\widehat{T}}(\omega)} \| \vartheta(t, \omega) \|^2 \,\mathrm{d}t =
\infty,\qquad  0 \le T < \mathcal{\widehat{T}}(\omega)
\end{equation}
for every $ \omega\in\{ \mathcal{\widehat{T}} < \infty\} $, and $
\int_0^T \| \vartheta(t, \omega) \|^2 \,dt < \infty$, $ \forall T
\in[0, \infty) $ for every $ \omega\in\{ \mathcal{\widehat{T}} =
\infty\} $. In conjunction with (\ref{3.2.G}), we obtain from (\ref
{Ttau111}) that
\[
\int_0^T \| \widetilde{\vartheta} (t, \omega) \|^2 \,\mathrm{d} t <
\int
_0^{\mathcal{\widehat{T}}(\omega)} \| \widetilde{\vartheta} (t,
\omega)
\|^2 \,\mathrm{d}t = \infty,\qquad  0 \le T < \mathcal{\widehat{T}}(\omega)
\]
holds for every $ \omega\in\{ \mathcal{\widehat{T}} < \infty\} $,
and that $ \int_0^T \| \widetilde{\vartheta} (t, \omega) \|^2
\,\mathrm{d} t < \infty$, $ \forall T \in[0, \infty) $ holds for every $
\omega\in\{ \mathcal{\widehat{T}} = \infty\} $.

We deduce for the stopping rules of (\ref{Ttau1}) and (\ref{Ttaui1})
the identification $ \mathcal{\widehat{T}}(\omega) = \mathcal
{T}(\omega
) $.

\subsection{The same thread, in reverse}
\label{subsec7.3}

Let us fix now a stochastic system $ \mathcal{N}\in\mathfrak{N}
(\mathbf{ x}) $ as in Section \ref{subsec7.2}, pick a progressively
measurable functional $ \alpha\dvtx[0, \infty) \times\Omega
\rightarrow
\mathbb{S}^n $ with $ \alpha(t, \omega) \in\mathcal{A} ( \omega
(t) )
$ for all $ (t,\omega) \in[0, \infty) \times\Omega$ and select a
progressively measurable functional $ \vartheta\dvtx[0, \infty)
\times
\Omega\rightarrow\mathbb{R}^n $ as in (\ref{ThetaAlpha}). For
\textit{this} $ \vartheta(\cdot, \cdot) $ and \textit{this} $
\mathcal{N}\in
\mathfrak{N} (\mathbf{ x}) $, we define $ \widetilde{\vartheta}
(\cdot
, \cdot) $ by (\ref{Ttau101}) as well as
%
\begin{eqnarray}
\label{LAMDA}
\Lambda(t) = \exp\biggl\{ \int_0^t ( \widetilde{\vartheta} (s, \mathfrak{
X} ) )' \,\mathrm{d} \widetilde{W}(s) - { 1 \over2} \int_0^t
\| \widetilde{\vartheta} (s, \mathfrak{ X} ) \|^2 \,\mathrm{d} s \biggr\}
\nonumber\\[-8pt]\\[-8pt]
\eqntext{\mbox{for } 0 \le t < \mathcal{T}}
\end{eqnarray}
as in (\ref{Lambda}), and set
%
\begin{equation}
\label{LAMDA_Too}
\Lambda(\mathcal{T} +u) =0 \qquad \mbox{for }
u \ge0 \mbox{ on } \{ \mathcal{T} < \infty\}
\end{equation}
in the notation of (\ref{Ttau1}). The resulting process
$ \Lambda(\cdot) $ is a local martingale and a supermartingale under
$ \mathbb{Q} $, and we have $
\mathcal{ T} (\mathfrak{ X}) = \inf\{ t \ge0 | \Lambda(t)=0\} $, $
\mathbb{Q}$-a.e.

We introduce also the sequence of $ \mathbb{F}$-stopping rules
\[
S_n (\omega) := \inf\biggl\{ t \ge0 \Big| \int_0^t
\| \widetilde{\vartheta} (s, \omega) \|^2 \,\mathrm{d}s \ge2 \log
n \biggr\} ,\qquad  n \in\mathbb{N} ,
\]
which satisfy
$
\lim_{n \rightarrow\infty} \uparrow S_n (\omega) =
\mathcal{T} (\omega) $ and $ \exp\{ {1 \over2 }
\int_0^{S_n (\omega)} \| \widetilde{\vartheta} (t, \omega)\|^2
\,\mathrm{d} t \} \le n$ $( \forall n \in
\mathbb{N}) $,
for every $ \omega\in\Omega$. From Novikov's theorem [e.g., \citet{KarShr91}, page 198], $
\Lambda( \cdot\wedge S_n) $ is a uniformly integrable $
\mathbb{Q}$-martingale; in particular, $ \mathbb{E}^\mathbb{Q} (
\Lambda( S_n)) =1 $ holds for every $ n \in\mathbb{N} $.
Thus, the recipe
\[
\mathbb{P}_n (A) := \mathbb{E}^\mathbb{Q} [ \Lambda( S_n) \cdot
\mathbf{ 1}_A ] ,\qquad  A \in\mathcal{F}(S_n)
\]
defines a consistent sequence, or ``tower,'' of probability measures
$ \{ \mathbb{P}_n \}_{n \in\mathbb{N}} $ on
$ (\Omega, \mathcal{F}) $. Appealing to the results in \citet{Par67}, pages 140--143 [see also the Appendix of F\"{o}llmer
(\citeyear{Fol72})], we
deduce the existence of a~probability measure $ \mathbb{P} $ on $
(\Omega, \mathcal{F}) $ such that
%
\begin{equation}
\label{tower}
\qquad \mathbb{P} (A) = \mathbb{P}_n (A) = \mathbb{E}^\mathbb{Q}
[ \Lambda( S_n) \cdot\mathbf{ 1}_A ] \quad \mbox{holds for every } A \in
\mathcal{F}(S_n) , n \in\mathbb{N}.\hspace*{-25pt}
\end{equation}
(Here again, the special structure imposed in this section on the
filtered measurable space $ (\Omega, \mathcal{F}), \mathbb{F} =
\{ \mathcal{F} (t) \}_{0 \le t < \infty} $ is indispensable.)
Therefore, for every $ T \in(0, \infty) $ we have
\[
\mathbb{P} (S_n >T) = \mathbb{E}^\mathbb{Q} \bigl[ \Lambda( S_n) \cdot
\mathbf{ 1}_{\{S_n >T\}} \bigr] = \mathbb{E}^\mathbb{Q} \bigl[ \Lambda( T)
\cdot
\mathbf{ 1}_{\{S_n >T\}} \bigr]
\]
by optional sampling, whereas monotone convergence leads to
%
\begin{equation}
\label{Wald}
\mathbb{P} (\mathcal{T} >T) = \mathbb{E}^\mathbb{Q} \bigl[ \Lambda( T)
\mathbf{ 1}_{\{\mathcal{T} >T\}} \bigr] .
\end{equation}
The following result echoes similar themes in \citet{CheFilYor05}.

\begin{lemma}
\label{Lemma_1}
The process $ \Lambda(\cdot) $ of (\ref{LAMDA}), (\ref{LAMDA_Too}) is
a $ \mathbb{Q}$-martingale, if and only if we have
\[
\mathbb{P}
(\mathcal{T} < \infty)=0
\]
[i.e., if and only if the process $
\mathfrak{ X}(\cdot) $ never hits the boundary of the orthant $
(0, \infty)^n $, $ \mathbb{P}$-a.s.].
\end{lemma}

\begin{pf}
If $ \mathbb{P} (\mathcal{T} < \infty)=0 $
holds, the nonnegativity of $ \Lambda(\cdot) $ and (\ref{Wald}) give
\[
1 = \mathbb{P} (\mathcal{T} >T) = \mathbb{E}^\mathbb{Q} \bigl[ \Lambda( T)
\mathbf{ 1}_{\{\mathcal{T} >T\}} \bigr] \le\mathbb{E}^\mathbb{Q} [
\Lambda( T) ] \qquad  \forall T \in(0, \infty) .\vadjust{\goodbreak}
\]
But $ \Lambda(\cdot) $ is a $ \mathbb{Q}$-supermartingale, so the
reverse inequality $ \mathbb{E}^\mathbb{Q} [ \Lambda( T) ] \le
\Lambda
(0) =1 $ also holds. We conclude that $ \mathbb{E}^\mathbb{Q} [
\Lambda
( T) ]=1 $ holds for all $ T \in(0, \infty) $, so $ \Lambda(\cdot) $
is a $ \mathbb{Q}$-martingale.

If, on the other hand, $ \Lambda(\cdot) $ is a $ \mathbb
{Q}$-martingale, then $ \mathbb{E}^\mathbb{Q} [ \Lambda( T) ]=1 $ and~%
(\ref{Wald}) give
\begin{eqnarray*}
\mathbb{P} ( \mathcal{T} \le T) &=& \mathbb{E}^\mathbb{Q} ( \Lambda( T)
) - \mathbb{E}^\mathbb{Q} \bigl( \Lambda(T) \mathbf{ 1}_{ \{ \mathcal
{T} >
T \}} \bigr) = \mathbb{E}^\mathbb{Q} \bigl( \Lambda(T) \mathbf{ 1}_{ \{
\mathcal
{T} \le T \}} \bigr) \\
&=& \mathbb{E}^\mathbb{Q} \bigl( \Lambda(\mathcal{T})
\mathbf
{ 1}_{ \{ \mathcal{T} \le T \}} \bigr)=0
\end{eqnarray*}
for every $ T \in[0, \infty) $, from Optional Sampling and the fact
that $ \Lambda(\mathcal{T}) = 0 $ holds $
\mathbb{Q} $-a.e. on $\{ \mathcal{T} <\infty\}$. We conclude $
\mathbb
{P} ( \mathcal{T} < \infty)=0 $; in conjunction with the identification
$ \mathcal{T} \equiv\widehat{\mathcal{T}} $ and (\ref{Ttaui1}), this
means that the co\"{o}rdinate mapping process $ \mathfrak{ X}(\cdot) $
never reaches the boundary of (i.e., takes values in) the strictly
positive orthant $(0, \infty)^n $, $ \mathbb{P}$-a.e.
\end{pf}

When the conditions of Lemma \ref{Lemma_1} prevail, the process
\[
W(\cdot) = \widetilde{W} (\cdot) - \int_0^\cdot\widetilde
{\vartheta}
( t, \mathfrak{ X} ) \,\mathrm{d} t
\]
is Brownian motion under the probability measure $ \mathbb{P} \equiv
\mathbb{P}^\mathcal{M} $ introduced\break in~(\ref{tower}). This measure
satisfies the equations of (\ref{Foell}), whereas the process~$
\mathfrak{ X}(\cdot) $ solves $ \mathbb{P}^\mathcal{M}$-a.s. the system
\[
d X_i(t) = X_i(t)
\sum_{\nu=1}^n \sigma_{i \nu}(t, \mathfrak{ X} ) [ \mathrm{d}
W_\nu(t) +
\vartheta_\nu(t, \mathfrak{ X}) \,\mathrm{d} t ] ,\qquad  X_i(0)=x_i>0
\]
for $ i=1,\ldots,n $, as in (\ref{1.1}).
It is then not hard to check that $ L(\cdot) $ defined by (\ref{expo})
satisfies $ \mathbb{P}^\mathcal{M}$-a.s. the identity $L(\cdot) X
(\cdot
) =( x_1 + \cdots+ x_n)/ \Lambda(\cdot) $ in accordance with (\ref{Lambda}).

We formalize these considerations as follows.

\begin{ass}\label{assC}
Suppose that the collection of sets $ \mathbb
{K} $ in (\ref{1.b}) satisfies (\ref{3.2.G}), (\ref{ThetaTrace}) and
that for any given progressively measurable functional $ \alpha\dvtx[0,
\infty) \times\Omega\rightarrow\mathbb{S}^n $ which satisfies
(\ref
{Alpha}) and
%
\begin{equation}
\label{AA1}
\int_0^T \operatorname{Tr} ( \alpha(t, \mathfrak{w}) ) \,\mathrm{d} t <
\infty\qquad \mbox{for all } (T, \mathfrak{ w}) \in[0, \infty)
\times\mathfrak{W} ,
\end{equation}
there exists a progressively measurable functional $ \vartheta\dvtx[0,
\infty) \times\Omega\rightarrow\mathbb{R}^n $ that satisfies the
condition (\ref{ThetaAlpha}), thus also by virtue of (\ref{ThetaTrace})
%
\begin{equation}
\label{AA2}
\int_0^T \| \vartheta(t, \mathfrak{ w}) \|^2 \,\mathrm{d} t < \infty
\qquad \mbox{for all } (T, \mathfrak{ w}) \in[0, \infty) \times
\mathfrak{W} .
\end{equation}
\end{ass}

The analysis of this subsection shows that, under Assumption \ref{assC} and
starting with any initial configuration $ \mathbf{ x} = (x_1, \ldots,
x_n )' \in(0, \infty)^n $ and with an arbitrary ``auxiliary''
admissible stochastic system $ \mathcal{N}= ( (\Omega, \mathcal{F})$,
$\mathbb{F}, \mathbb{Q} , \mathfrak{X}(\cdot), \widetilde{W}(\cdot
) ) $
in $ \mathfrak{ N} (\mathbf{ x}) $ as in Section \ref{subsec7.2},
the process $ \Lambda(\cdot) $ of (\ref{LAMDA_Too}) is a $ \mathbb
{Q}$-martingale, and we can construct a ``primal'' admissible system $
\mathcal{ M} \in\mathfrak{ M} (\mathbf{ x}) $ as in Section \ref
{subsec2.1} [i.e., with the canonical process $ \mathfrak{X} (\cdot) $
taking values in $ (0, \infty)^n $, $ \mathbb{P}^\mathcal{M}$-a.s.],
for which (\ref{Foell}) holds, and we have $ \mathbb{P}^\mathcal{M}
\ll
\mathbb{Q} $ as in (\ref{bigtriangleleft}). We deduce
\[
Q (T, \mathbf{ x}) = \sup_{ \mathcal{N} \in\mathfrak{ N} (\mathbf{
x})} \mathbb{ Q}^{ {\mathcal{N}}} ( \mathcal{T} >T ) \le\sup
_{\mathcal
{ M} \in\mathfrak{ M} (\mathbf{ x})} \biggl( \frac{ \mathbb{E}^{ \mathbb
{P}^\mathcal{M}}[ L(T) X(T) ] }{ x_1 + \cdots+ x_n } \biggr)= \Phi(T,
\mathbf
{ x}) .
\]
The reverse inequality $ Q (T,\mathbf{x}) \ge\Phi(T,\mathbf{x}) $
was established in (\ref{3.7.f}). This way, for every function $ U
\dvtx
[0,\infty) \times(0, \infty)^n \rightarrow(0, \infty) $ in the
collection $ \mathcal{U} $, we can strengthen (\ref{2.8}) to
%
\begin{eqnarray}
\label{3.8}
U (T,\mathbf{x}) &\ge&\mathfrak{u} (T,\mathbf{x}) \ge\Phi
(T,\mathbf{x}) = Q (T,\mathbf{x}) \nonumber\\[-8pt]\\[-8pt]
\eqntext{\forall
(T,\mathbf{x}) \in(0,\infty) \times(0,\infty)^n .}
\end{eqnarray}
We have established the following result.

\begin{prop}
\label{Proposition_4}
Recall the functions $ \mathfrak{u}(\cdot, \cdot) $, $ \Phi(\cdot,
\cdot) $ and $ Q (\cdot, \cdot) $ defined on $ (0,\infty) \times
(0,\infty)^n $ by (\ref{1.7}), (\ref{1.8}) and (\ref{3.7.f}),
respectively, and impose Assumption~\ref{assC}. Then (\ref{3.8}) holds for every
function $ U (\cdot, \cdot) \in\mathcal{U} $.
\end{prop}

\begin{remark}
\label{Remark_4}
Here is a situation where Assumption \ref{assC} prevails:
Suppose that (\ref{3.2.G}) holds and that, for every $ \mathbf{ z}
\in
(0, \infty)^n $ and $ a \in\mathcal{A} (\mathbf{ z}) $, we have
$ ( \theta, a) \in\mathcal{K} (\mathbf{ z}) $ for $ \theta$ given by
$ \theta_\nu= \sum_{j =1}^n s_{j \nu} $, $ \nu=1, \ldots, n $ and
$ s
s' = a $. Then for any progressively measurable $ \alpha\dvtx[0,
\infty)
\times\Omega\rightarrow\mathbb{S}^n $ that satisfies (\ref{AA1}) we
select the progressively measurable functional $ \vartheta\dvtx[0,
\infty
) \times\Omega\rightarrow\mathbb{R}^n $ via $ \vartheta_\nu(t,
\omega) = \sum_{j =1}^n \sigma_{j \nu} (t, \omega) $, $ \nu=1,
\ldots,
n $. This choice induces
\[
\widetilde{\vartheta}_\nu(t, \omega) = \sum_{i=1}^n \biggl( 1 - { \omega_i
(t) \over\omega_1 (t) + \cdots+ \omega_n (t) } \biggr) \sigma_{i \nu} (t,
\omega) ,\qquad  \nu= 1, \ldots, n
\]
which obeys $ \int_0^T \| \widetilde{\vartheta} (t, \mathfrak{ w})
\|^2\,d t < \infty$ as in (\ref{AA2}) for all $ (T, \mathfrak{ w}) \in[0,
\infty) \times\mathfrak{W} $; the process $ \Lambda(\cdot) $ of
(\ref{LAMDA}) and (\ref{LAMDA_Too}) is a $ \mathbb{Q}$-martingale, whereas
(\ref{1.1}) becomes
\[
\mathrm{d} X_i(t) = X_i(t) \Biggl[
\sum_{\nu=1}^n \sigma_{i \nu}(t, \mathfrak{ X} ) \,\mathrm{d} W_\nu
(t) +
\Biggl( \sum_{j=1}^n \alpha_{i j}(t, \mathfrak{ X} ) \Biggr)
\,\mathrm{d} t \Biggr] ,\quad  X_i(0)=x_i>0
\]
for $ i=1, \ldots, n $. The condition (\ref{AA1}) guarantees now that $
\mathfrak{ X} (\cdot) $ takes values in $ (0, \infty)^n$, $ \mathbb
{P}^\mathcal{M}$-a.s. in the resulting primal admissible system $
\mathcal{ M} \in\mathfrak{ M} (\mathbf{ x}) $.
\end{remark}

\section{Dynamic programming}\label{sec8}

The quantity $ Q
(T,\mathbf{x}) $ defined in (\ref{3.7.f}) is the value of a stochastic
control problem: namely, the \textit{maximal}
``\textit{containment}'' \textit{probability}, over all measures $ \mathbb{Q}^\mathcal{N}
$ with $ \mathcal{N} \in\mathfrak{ N} (\mathbf{ x}) $, that the
process $ \mathfrak{X} (\cdot) $ with dynamics (\ref{aux}),
initial configuration $ \mathfrak{X} (0) =\mathbf{x} \in(0,
\infty)^n $, and controlled through the choice of progressively
measurable functional $ \alpha(\cdot, \cdot) $ as in
(\ref{const1}), (\ref{const2}), does \textit{not} hit the boundary of the
positive orthant by time $ T$.

Let us suppose that the resulting function $ Q (\cdot, \cdot) $ is
continuous on $ (0, \infty) \times(0, \infty)^n $. Then it can be
checked [as in Lions (\citeyear{Lio84}), Lemma II.1 and \citet{Lio83N1}, Theorem
II.4] that it satisfies as well the following dynamic programming
principle: \textit{for every initial configuration $ \mathbf{ x} \in(0,
\infty)^n$, the process}
%
\begin{eqnarray}
\label{DP}
Q \bigl(T-t , \mathfrak{ X}( t ) \bigr) \mathbf{ 1}_{\{ \mathcal{T} >t\}} ,\nonumber\\[-8pt]\\[-8pt]
\eqntext{0\le t \le T \mbox{ \textit{is a} } \mathbb{Q}^{ \mathcal{N}}\mbox{-\textit{supermartingale}, } \forall\mathcal{N} \in\mathfrak{
N}(\mathbf{ x}) .}
\end{eqnarray}
[See \citet{ElKHuuJea87}, \citet{HauLep90},
\citet{FleSon93} and \citet{Kry80}
for results in a similar vein.] Equivalently, the process $ L(t) X(t) Q
(T-t , \mathfrak{ X}( t ) ) $, $ 0 \le t \le T $ is a $ \mathbb
{P}^{ \mathcal{M}}$-supermartingale, for every $ \mathcal{M} \in
\mathfrak{M}(\mathbf{ x}) $.

Consider now an arbitrary continuous function $ \breve{U}\dvtx[0,
\infty)
\times(0, \infty)^n \rightarrow[0, \infty) $ which satisfies $
\breve
{U}(0, \cdot) \equiv1 $ on $(0, \infty)^n $, and is such that, for
every $ \mathcal{N} \in\mathfrak{ N}(\mathbf{ x}) $ and $ \mathbf{ x}
\in(0, \infty)^n $, the process $ \breve{U} (T-t , \mathfrak{ X}( t )
) \mathbf{ 1}_{\{ \mathcal{T} >t\}} , 0 \le t \le T $ is a~$ \mathbb
{Q}^{ \mathcal{N}}$-super\-martingale. We shall denote by $ \mathcal
{\breve{U}} $ the collection of all such functions and note that $
\mathcal{U} \subseteq\mathcal{\breve{U}} $ and $ Q \in\mathcal
{\breve
{U}} $.
From optional sampling we have then for every $ \mathcal{N} \in
\mathfrak{ N}(\mathbf{ x}) $ the comparisons
\[
\breve{U}(T, \mathbf{ x}) \ge\mathbb{E}^{ \mathbb{Q}^{ \mathcal{N}}}
\bigl[ \breve{U} ( 0, \mathfrak{ X} (T ) ) \mathbf{ 1}_{\{ \mathcal{T}
>T \}
} \bigr] = \mathbb{Q}^{ \mathcal{N}} ( \mathcal{T} >T ) ,
\]
thus also $ \breve{U}(T, \mathbf{ x}) \ge Q(T, \mathbf{ x}) $, $
\forall(T, \mathbf{ x}) \in(0, \infty) \times(0, \infty)^n $.

In other words, \textit{the function $ Q (\cdot,\cdot) $ defined in}
(\ref{3.7.f}) \textit{is the smallest element of the collection} $ \mathcal{\breve
{U}} $. It is also clear from this line of reasoning that $ \mathcal
{N}_o \in\mathfrak{N} (\mathbf{ x}) $ attains the supremum $ Q (T,
\mathbf{ x}) =\sup_{ \mathcal{N} \in\mathfrak{N} (\mathbf{ x})}
\mathbb{
Q}^{ \mathcal{N}} ( \mathcal{T} >T ) $ in (\ref{3.7.f}), if and only if
the process $ Q (T- t , \mathfrak{ X}( t ) ) \mathbf{ 1}_{\{ \mathcal
{T} >t \}} , 0 \le t \le T $ is a $ \mathbb{Q}^{ \mathcal{N}_o}$-martingale.

\begin{thm}
\label{Theorem_1}
Suppose that Assumption \ref{assC} and conditions (\ref{measur})--(\ref
{UppSemiCont}) hold and that the function $ Q (\cdot, \cdot) $ of
(\ref
{3.7.f}) is continuous on $ (0, \infty) \times(0, \infty)^n $. Then
the infimum in (\ref{1.7}) is attained, and
%
\begin{equation}
\label{3.8.a}
\mathfrak{u} (T,\mathbf{x}) = \Phi(T,\mathbf{x}) = Q
(T,\mathbf{x}) \qquad  \forall(T,\mathbf{x}) \in(0,\infty) \times
(0,\infty
)^n .
\end{equation}
\end{thm}

\begin{pf}
Consider an arbitrary function $ \breve{U} (\cdot,
\cdot
) $ in the collection $ \mathcal{\breve{U}} $ just defined and fix an
arbitrary pair $ (T,\mathbf{x}) \in(0,\infty)
\times(0,\infty)^n $; then for every $ \varepsilon>0 $, consider a
mollification $ U_\varepsilon(\cdot, \cdot) \in\mathcal{U} $ of the
function $ \breve{U} (\cdot, \cdot) $ with $ U_\varepsilon
(T,\mathbf
{x}) \le\breve{U} (T,\mathbf{x}) + \varepsilon$.

Proposition \ref{Proposition_2} gives then $ \mathfrak{ u} (T,\mathbf
{x}) \le\breve{U} (T,\mathbf{x}) + \varepsilon$. Since $
\varepsilon
>0 $ is arbitrary, this shows that $ \mathfrak{ u} (T,\mathbf{x}) $ is
dominated by $ Q (T,\mathbf{x}) $, the infimum of $ \breve{U}
(T,\mathbf
{x}) $ over all functions $ \breve{U} (\cdot, \cdot) \in\mathcal
{\breve{U}} $. But the reverse inequality $ \mathfrak{u} (T,\mathbf{x})
\ge Q (T,\mathbf{x}) $ holds on the strength of (\ref{3.8}), so (\ref
{3.8.a}) follows.
\end{pf}

\subsection{The HJB equation}
\label{subsec8.1}

Under the conditions of Theorem \ref{Theorem_1}, the arbitrage function
$ \mathfrak{ u}(\cdot, \cdot) $ is equal to the function $ Q(\cdot,
\cdot) $ of (\ref{3.7.f}) and is continuous on $ (0, \infty) \times(0,
\infty)^n $. Thanks to the dynamic programming principle of~(\ref{DP}),
it is also a viscosity solution of the
\textit{Hamilton--Jacobi--Bellman} (\textit{HJB}) \textit{equation}
%
\begin{equation}
\label{3.9}
\qquad\frac{\partial U}{\partial\tau} (\tau, \mathbf{z}) = \sup_{a
\in\mathcal{A} (\mathbf{z})} \sum_{i=1}^n \sum_{j=1}^n
z_i z_j a_{ij} \biggl( \frac{1}{ 2 } D^2_{ij} U(\tau,
\mathbf{z})+ \frac{D_{i} U(\tau, \mathbf{z})}{z_1 + \cdots+z_n}
\biggr)
\end{equation}
on $ (0, \infty) \times(0, \infty)^n $ [cf. Lions (\citeyear{Lio84}), Theorem
III.1 or \citet{Lio83N2}, Theorem~I.1].

If in addition to being continuous, as we assumed in Theorem \ref
{Theorem_1}, the function $ Q(\cdot, \cdot) $ of (\ref{3.7.f}) is of
class $ \mathcal{C}^{1,2} $ locally on $ (0, \infty) \times(0,
\infty)^n
$, then the arbitrage function $ \mathfrak{ u}(\cdot, \cdot) $ is not
only a viscosity solution but actually a~classical solution of the HJB
equation (\ref{3.9}). This is the case, for instance, under the
combined conditions of Theorem \ref{Theorem_1} and Proposition \ref
{Proposition_3}; then the arbitrage function $ \mathfrak{ u}(\cdot, \cdot)
$ also satisfies on the domain $ (0, \infty) \times(0, \infty)^n $ the
\textit{linear} parabolic equation
%
\begin{equation}
\label{7.25}
\frac{\partial U}{\partial\tau} (\tau, \mathbf{z}) = \sum_{i=1}^n
\sum_{j=1}^n z_i z_j \mathbf{a}_{ij} (\mathbf{ z}) \biggl( \frac{1}{ 2 }
D^2_{ij} U (\tau,\mathbf{z})+ \frac{D_{i} U (\tau, \mathbf{z})}{z_1+
\cdots+z_n} \biggr)
\end{equation}
with $ \mathbf{a}\dvtx(0, \infty)^n \rightarrow\mathbb{S}^n $ as in
Assumption \ref{assA} or Remark \ref{Remark_2}, in addition to the initial condition
%
\begin{equation}
\label{3.9.a}
U (0, \cdot) \equiv1 \qquad \mbox{on } (0, \infty)^n .
\end{equation}
In particular, the arbitrage function $ \mathfrak{ u}(\cdot, \cdot) $
satisfies, in this case, the requirement (\ref{D.11}) and belongs to the
class $ \mathcal{U} $ of Section \ref{sec5}.

Recalling Propositions \ref{Proposition_1}--\ref{Proposition_3} and
Theorem \ref{Theorem_1}, we summarize the above discussion as follows.

\begin{thm}
\label{Theorem_2} Suppose that conditions (\ref{measur})--(\ref
{UppSemiCont}), (\ref{3.2.GG}), (\ref{3.2}) and Assumptions \ref{assA}, \ref{assB} and \ref{assC}
are in force.

Then the arbitrage function $ \mathfrak{ u}(\cdot, \cdot) $ is the
smallest element of the class $ \mathcal{U} $, as well as a classical
solution of both the HJB equation (\ref{3.9}) and of the linear
parabolic equation (\ref{7.25}), subject to (\ref{3.9.a}). Furthermore
(\ref{3.8.a}) holds, the infimum in (\ref{1.7}) is attained, and the
Markovian investment rule $ \pi^U (\cdot, \cdot) $ in (\ref{2.11.a})
with $ U \equiv\mathfrak{ u} $ satisfies (\ref{2.11.b}) for every
admissible system $ \mathcal{M} \in\mathfrak{M} (\mathbf{ x}) $.
\end{thm}

\begin{remark}
\label{Remark_5}
We note that Theorem \ref{Theorem_2} is in
agreement with general regularity theory for fully nonlinear parabolic
equations, as in Lions (\citeyear{Lio83N3}), Theorem II.4 (see also \citet{Kry87},
Section 6.5; Krylov (\citeyear{Kry89N2}); Wang (\citeyear{WAN92}, \citeyear{Wan92N2}, \citeyear{Wan92N3}), Theorems II.3.2 and III.2; or
\citet{Lie96}, Chapter XIV).

As we mentioned already, Assumptions \ref{assA} and \ref{assB} can be replaced in Theorem
\ref{Theorem_2} by the conditions of Remark \ref{Remark_2}. We
conjecture that the conclusions of Theorem \ref{Theorem_2} should hold
under even weaker assumptions but leave this issue for future research.

We also remark that the function $ V (t, \mathbf{ z}) := (z_1 + \cdots
+ z_n) \mathfrak{ u} (t, \mathbf{ z}) $, $ (\tau, z) \in(0, \infty)
\times(0, \infty)^n $ satisfies an HJB-type equation simpler than
(\ref
{3.9}), namely, the \textit{Pucci maximal equation},
%
\begin{equation}
\label{3.10}
\frac{\partial V }{\partial\tau} (\tau, \mathbf{z}) = \frac{1}{ 2 }
\sup_{a
\in\mathcal{A} (\mathbf{z})} \sum_{i=1}^n \sum_{j=1}^n
z_i z_j a_{ij} D^2_{ij} V (\tau,\mathbf{z}) ,
\end{equation}
along with the initial condition
$ V (0, \mathbf{z}) = z_1 + \cdots+ z_n $.
In the setting of Theorem~\ref{Theorem_2}, equation (\ref{3.10})
reduces to
\[
\frac{\partial V }{\partial\tau} (\tau, \mathbf{z}) = \frac{1}{ 2 }
\sum_{i=1}^n \sum_{j=1}^n
z_i z_j \mathbf{a}_{ij} (\mathbf{ z}) D^2_{ij} V (\tau,\mathbf{z}) .
\]
\end{remark}

\subsection{An example}
\label{subsec8.2}

Let us go back to the volatility-stabilized model introduced in
\citet{FERKAR05}, but now with some ``Knightian''
uncertainty regarding its volatility structure
\[
1 \le\alpha_{ii} (t) \mu_i (t) \le1+ \delta,\qquad  0 \le t < \infty
\]
for some given $ \delta\ge0 $. The case $\delta=0$ corresponds to
the variance structure of the model studied in Fernholz and Karatzas (\citeyear{FERKAR05}, \citeyear{FERKAR}).

More specifically let us assume that, for any given $ \mathbf{y} \in
\mathfrak{S}_n $, the compact, convex subset $ \mathcal{A} (\mathbf{y})
$ of $ \mathbb{S}^n $ in (\ref{E.3.a})
consists of all matrices $ a = \{ a_{ij} \}_{1 \le i,j \le n} $
with $ a_{ij} =0 $ for $ j \neq i $ and
\[
y_i a_{ii} = \eta^2 (y_1 + \cdots+ y_n ) ;\qquad  i=1, \ldots, n , 1 \le
\eta
\le1 + \delta.
\]
The sets of (\ref{1.b}) are given as
\[
\mathcal{K} (\mathbf{ y}) = \bigl\{ (\theta, a) | a \in\mathcal{A}
(\mathbf
{ y}), \theta= \bigl( \zeta\sqrt{ a_{11} }, \ldots, \zeta\sqrt{ a_{nn} }
\bigr)' \mbox{ with } \zeta\in\bigl[ \sqrt{C_1 }, \sqrt{C_2 } \bigr] \bigr\}
\]
for some given constants $ C_1 \in(0, 1] $, $ C_2 \in(1, \infty) $;
these choices satisfy (\ref{3.2.G}), (\ref{ThetaTrace}).

Condition (\ref{3.2}) is satisfied in this case automatically
(in fact, with $ \lambda\equiv1$), as are 
(\ref{3.2.GG}) and (\ref{D.4}):
it suffices to take
%
\begin{equation}
\label{VolStab}
\mathbf{a}_{ii} (\mathbf{ z}) = (z_1 + \cdots+ z_n) / z_i ,\qquad  i = 1,\ldots, n
\end{equation}
and $ H(\mathbf{ z}) = \sum_{i=1}^n \log z_i $, which induces $ \theb_i
(z) = \sqrt{\mathbf{a}_{ii} (\mathbf{ z}) } $ in (\ref{D.2}). These
functions are all locally bounded and locally Lipschitz continuous on $
(0, \infty)^n $.

The HJB equation (\ref{3.9}) satisfied by the arbitrage function $
\mathfrak{ u}(\cdot, \cdot) $ becomes
\[
\frac{\partial U}{\partial\tau} (\tau, \mathbf{z}) = \sup_{1 \le
\eta
\le1+\delta} \Biggl[ \eta^2 \Biggl\{ \frac{1}{ 2 } \sum_{i=1}^n (z_1 + \cdots+
z_n) z_i D^2_{ii} U(\tau,
\mathbf{z})+\sum_{i=1}^n z_i D_{i} U(\tau, \mathbf{z}) \Biggr\} \Biggr] ,
\]
and reduces to the \textit{linear} parabolic equation
%
\begin{equation}
\label{3.11}
\frac{\partial U}{\partial\tau} (\tau, \mathbf{z}) =
\frac{1}{ 2 } \sum_{i=1}^n (z_1 + \cdots+ z_n) z_i
D^2_{ii} U (\tau,
\mathbf{z})+\sum_{i=1}^n z_i D_{i} U (\tau, \mathbf{z})
\end{equation}
of (\ref{7.25}) for the choice of variances in (\ref{VolStab}). The
reason for this reduction is that the expression on the left-hand side
of (\ref{3.11}) is negative, so we have considerable simplification in
this case.

\begin{remark}
In this example, the arbitrage function $ \mathfrak{ u}(\cdot, \cdot)
$ can be represented as
\[
\mathfrak{ u} (T, \mathbf{ z}) = { z_1 \cdots z_n \over z_1+ \cdots+
z_n } \mathbb{E} \biggl[ { X_1 (T) + \cdots+ X_n (T) \over X_1 (T) \cdots
X_n (T)} \biggr]
\]
in terms of the components of the $ (0, \infty)^n$-valued
capitalization process $ \mathfrak{ X} (\cdot) = ( X_1 (\cdot),
\ldots
, X_n (\cdot) )' $. These are now time-changed versions $ X_i (\cdot)
=\break
\Psi_i ( A (\cdot) ) $, $ i=1, \ldots, n $ of the independent
squared-Bessel processes
\[
\mathrm{d} \Psi_i (u) = 4 u \,\mathrm{d} u + 2 \sqrt{ \Psi_i (u) }
\,\mathrm{d} \beta_i (u) ,\qquad  \Psi_i (0) = z_i ,
\]
run with a time change $ A( t) = (1 / 4) \int_0^t (X_1 (s) + \cdots+
X_n (s) ) \,\mathrm{d}s $ common for all components, and with $ \beta_1
(\cdot), \ldots, \beta_n (\cdot) $ independent standard Brownian
motions [see Fernholz and Karatzas (\citeyear{FERKAR05}, \citeyear{FERKAR}), \citet{Goi09} and Pall
(\citeyear{PAL}) for more details].
\end{remark}

\section{A stochastic game}
\label{sec9}

For any given investment rule $ \Pi\in\mathfrak{ P} $ and admissible
system $ \mathcal{M} \in\mathfrak{M} (\mathbf{ x}) $, let us
consider the quantity
%
\begin{equation}
\label{9.1}
\xib_{ \Pi, \mathcal{M}} (T,\mathbf{x}) : = \inf\bigl\{ r>0 \dvtx
\mathbb
{P}^\mathcal{M} \bigl( Z^{ r X(0), \Pi} ( T) \ge X(T) \bigr) = 1 \bigr\} .
\end{equation}
This measures, as a proportion of the initial total market
capitalization, the smallest initial capital that an investor who uses
the rule $ \Pi$ and operates within the market model $ \mathcal{M} $,
needs to set aside at time $ t=0 $ in order for his wealth to be able
to ``catch up with the market portfolio'' by time $ t = T $, with $
\mathbb{P}^\mathcal{M}$-probability one.

Our next result exhibits the arbitrage function $ \mathfrak{ u}(\cdot
, \cdot) $ of (\ref{1.7}) as the min--max value of a zero-sum stochastic
game between two players: the investor, who tries to select the rule $
\Pi\in\mathfrak{P} $ so as to make the quantity of (\ref{9.1}) as
small as possible and ``nature,'' or the goddess Tyche herself, who
tries to thwart him by choosing the admissible system or ``model'' $
\mathcal{M} \in\mathfrak{M} (\mathbf{ x}) $ to his detriment.

\begin{thm}
\label{Theorem_3}
Under the conditions of Theorem \ref{Theorem_1}, we have
%
\begin{equation}
\label{9.2}
\mathfrak{u} (T,\mathbf{x}) = \inf_{\Pi\in\mathfrak{ P}} \Bigl( \sup
_{\mathcal{M} \in\mathfrak{M} (\mathbf{ x})} \xib_{ \Pi, \mathcal{M}}
(T,\mathbf{x}) \Bigr) = \sup_{\mathcal{M} \in\mathfrak{M} (\mathbf{
x})} \Bigl( \inf
_{\Pi\in\mathfrak{ P}} \xib_{ \Pi, \mathcal{M}} (T,\mathbf{x}) \Bigr) .
\end{equation}
\end{thm}

\begin{pf}
For the quantities of (\ref{9.1}) and (\ref{1.9}) we claim
%
\begin{equation}
\label{9.3}
\xib_{ \Pi, \mathcal{M}} (T,\mathbf{x}) \ge\mathfrak{ u}_{
\mathcal
{M}} (T,\mathbf{x})  \qquad \forall(\Pi, \mathcal{M}) \in\mathfrak{ P}
\times\mathfrak{M} (\mathbf{ x}) .
\end{equation}
Indeed, if the set on the right-hand side of (\ref{9.1}) is empty, we
have\break $ \xib_{ \Pi, \mathcal{M}} (T,\mathbf{x})= \infty$ and nothing
to prove; if, on the other hand, this set is not empty, then for any of
its elements $ r >0 $ the process
$ L( \cdot) V^{ r X(0), \Pi} ( \cdot) $ is a $ \mathbb{P}^\mathcal
{M}$-supermartingale, and therefore (\ref{1.89}), that is, $ r \ge
\mathfrak{ u}_{ \mathcal{M}} (T,\mathbf{x}) $, still holds and~(\ref
{9.3}) follows again.

Taking the infimum with respect to $ \Pi\in\mathfrak{ P} $ on the
left-hand side of (\ref{9.3}), then the supremum of both sides with
respect to $ \mathcal{M} \in\mathfrak{M} (\mathbf{ x}) $, we obtain
%
\[
\underline{G}(T,\mathbf{x}) := \sup_{\mathcal{M} \in\mathfrak{M}
(\mathbf
{ x})} \Bigl( \inf_{\Pi\in\mathfrak{ P}} \xib_{ \Pi, \mathcal{M}}
(T,\mathbf{x}) \Bigr)
\ge\sup_{\mathcal{M} \in\mathfrak{M} (\mathbf{ x})} \mathfrak{ u}_{
\mathcal{M}} (T,\mathbf{x}) = \Phi(T,\mathbf{x})
\]
from (\ref{1.9}). The quantity $ \underline{G}(T,\mathbf{x}) $ is the
lower value of the stochastic game under consideration.

In order to complete the proof of (\ref{9.2}) it suffices, on the
strength of Theorem \ref{Theorem_1}, to show that the upper value
\[
\overline{G}(T,\mathbf{x}) := \inf_{\Pi\in\mathfrak{ P}} \Bigl( \sup
_{\mathcal{M} \in\mathfrak{M} (\mathbf{ x})} \xib_{ \Pi, \mathcal{M}}
(T,\mathbf{x}) \Bigr) \ge\underline{G}(T,\mathbf{x})
\]
of this game satisfies
%
\begin{equation}
\label{9.5}
\overline{G}(T,\mathbf{x}) \le\mathfrak{u} (T,\mathbf{x}) .
\end{equation}
To see this, we introduce for each given investment rule $ \Pi\in
\mathfrak{ P} $ the quantity
\[
\mathfrak{h}_{ \Pi} (T,\mathbf{x}) := \inf\bigl\{ r>0 \dvtx\mathbb
{P}^\mathcal
{M} \bigl( Z^{ r X(0), \Pi} ( T) \ge X(T) \bigr) = 1 ,
\forall\mathcal{M} \in\mathfrak{M} (\mathbf{x}) \bigr\} ;
\]
that is, the smallest proportion $ r >0 $ of the initial market capitalization
that allows an investor using the rule $ \Pi$ to be able to ``catch up
with the market portfolio'' by time $ t = T $ with $ \mathbb
{P}^\mathcal
{M}$-probability one, no matter which admissible system (model) $
\mathcal{M} $ might materialize. We have clearly
%
\begin{equation}
\label{9.7}
\mathfrak{h}_{ \Pi} (T,\mathbf{x}) \ge\mathfrak{u} (T,\mathbf{x})
\vee
\xib_{ \Pi, \mathcal{M}} (T,\mathbf{x}) \qquad  \forall(\Pi, \mathcal{M})
\in\mathfrak{ P} \times\mathfrak{M} (\mathbf{ x}) ,
\end{equation}
which leads to
%
\begin{equation}
\label{9.8}
\mathfrak{u} (T,\mathbf{x}) = \inf_{\Pi\in\mathfrak{ P} }
\mathfrak
{h}_{ \Pi} (T,\mathbf{x}) \ge
\inf_{\Pi\in\mathfrak{ P}} \Bigl( \sup_{\mathcal{M} \in\mathfrak{M}
(\mathbf
{ x})} \xib_{ \Pi, \mathcal{M}} (T,\mathbf{x}) \Bigr) = \overline
{G}(T,\mathbf{x})
\end{equation}
and proves (\ref{9.5}).
\end{pf}

\subsection{A least favorable model and the investor's best response}
\label{subsec9.1}

Let us place ourselves now in the context of Theorem \ref{Theorem_2}
and observe that Proposition \ref{Proposition_1}, along with
Proposition \ref{Proposition_2} and its Corollary, yields
%
\begin{equation}
\label{9.9}
\quad \xib_{ \Pi, \mathcal{M}_o} (T,\mathbf{x}) \ge\mathfrak{ u}_{
\mathcal
{M}_o} (T,\mathbf{x}) = \Phi(T,\mathbf{x}) = \xib_{ \Pi_o,
\mathcal
{M}_0} (T,\mathbf{x}) \qquad  \forall\Pi\in\mathfrak{ P} ,
\end{equation}
by virtue of (\ref{9.3}) for $ \mathcal{M} \equiv\mathcal{M}_o $ and
of (\ref{2.11.b}) for $ U \equiv\Phi$. Here $ \mathcal{M}_o $ is the
``least favorable admissible system'' that attains the supremum over $
\mathfrak{M} (\mathbf{ x}) $ in (\ref{1.8}), and $ \Pi_o \equiv\pi
^\Phi$ denotes the investment rule of (\ref{2.11.a}) with $ U \equiv
\Phi$.

In this setting, the investment rule $ \Pi_o \in\mathfrak{ P} $
attains the infimum\break $ \inf_{\Pi\in\mathfrak{ P} } \mathfrak{h}_{
\Pi}
(T,\mathbf{x}) = \mathfrak{u} (T,\mathbf{x}) $ in (\ref{9.8}), and we
obtain then
%
\begin{eqnarray}
\label{9.10}
\mathfrak{h}_{ \Pi_o} (T,\mathbf{x}) &=& \mathfrak{u} (T,\mathbf{x}) =
\Phi(T,\mathbf{x}) = \xib_{ \Pi_o, \mathcal{M}_0} (T,\mathbf{x})
\ge\xib_{ \Pi_o, \mathcal{M}} (T,\mathbf{x}) \nonumber\\[-8pt]\\[-8pt]
\eqntext{\forall
\mathcal{M} \in\mathfrak{M} (\mathbf{ x})}
\end{eqnarray}
on the strength of (\ref{9.7}). Putting (\ref{9.9}) and (\ref{9.10})
together we deduce
%
\begin{eqnarray}
\xib_{ \Pi, \mathcal{M}_o} (T,\mathbf{x}) \ge\mathfrak{ u}
(T,\mathbf
{x}) = \xib_{ \Pi_o, \mathcal{M}_0} (T,\mathbf{x}) \ge\xib_{ \Pi_o,
\mathcal{M}} (T,\mathbf{x}) \nonumber \\
\eqntext{\forall(\Pi, \mathcal{M}) \in
\mathfrak
{ P} \times\mathfrak{M} (\mathbf{ x}) ,}
\end{eqnarray}
the saddle property of the pair $ (\Pi_o, \mathcal{M}_o) \in
\mathfrak
{ P} \times\mathfrak{M} (\mathbf{ x}) $.

In particular, the investment rule $ \Pi_o \equiv\pi^\Phi$ of (\ref
{2.11.a}) with $ U \equiv\Phi$ is seen to be the investor's best
response to the least favorable admissible system $ \mathcal{M}_o \in
\mathfrak{M} (\mathbf{ x}) $ of Proposition \ref{Proposition_1}, and
vice-versa. In this sense the investor, once he has figured out a least
favorable admissible system $ \mathcal{M}_o $, can allow himself the
luxury to ``forget'' about model uncertainty and concentrate on finding
an investment rule $ \Pi_o \in\mathfrak{ P} $ that satisfies $ \xib_{
\Pi, \mathcal{M}_o} (T,\mathbf{x}) \ge\xib_{ \Pi_o, \mathcal{M}_o}
(T,\mathbf{x}) $, $ \forall\Pi\in\mathfrak{ P} $ as in (\ref{9.9}),
that is, on outperforming the market portfolio with the least initial
capital within the context of the least favorable model $ \mathcal
{M}_o $.

\section*{Acknowledgments}
We are grateful to Professors Jak\v{s}a Cvitani\'{c}, Nicolai Krylov,
Mete Soner, Nizar Touzi, Hans F\"{o}llmer, Erhan Bayraktar, Johan
Tysk, and to the two reviewers, for helpful advice on the subject
matter of this paper and for bringing relevant literature to our
attention. Numerous helpful discussions with Drs. Robert Fernholz,
Adrian Banner, Vasileios Papathanakos, Phi-Long Nguyen-Thanh, and
especially with Johannes Ruf, are gratefully acknowledged.



\printaddresses


\begin{thebibliography}{78}

\bibitem[\protect\citeauthoryear{Avellaneda, L{\'e}vy and
  Par{\'a}s}{1995}]{AVELEVPAR95}
\begin{barticle}[auto:STB|2011-03-03|12:04:44]
\bauthor{\bsnm{Avellaneda},~\bfnm{M.}\binits{M.}},
  \bauthor{\bsnm{L{\'e}vy},~\bfnm{M.~M.}\binits{M.~M.}} \AND
  \bauthor{\bsnm{Par{\'a}s},~\bfnm{A.}\binits{A.}}
(\byear{1995}).
\btitle{Pricing and hedging derivative securities in markets with uncertain
  volatility}.
\bjournal{Appl. Math. Finance}
\bvolume{2}
\bpages{73--88}.
\end{barticle}
\endbibitem

\bibitem[\protect\citeauthoryear{Bayraktar and Yao}{2011}]{BAYYAO}
\begin{bmisc}[auto:STB|2011-03-03|12:04:44]
\bauthor{\bsnm{Bayraktar},~\bfnm{E.}\binits{E.}} \AND
  \bauthor{\bsnm{Yao},~\bfnm{S.}\binits{S.}}
(\byear{2011}).
\bhowpublished{Optimal stopping for nonlinear expectations. \textit{Stochastic
  Process. Appl.} To appear}.
\end{bmisc}
\endbibitem

\bibitem[\protect\citeauthoryear{Bayraktar, Karatzas and Yao}{2011}]{BAYKARYAO}
\begin{bmisc}[auto:STB|2011-03-03|12:04:44]
\bauthor{\bsnm{Bayraktar},~\bfnm{E.}\binits{E.}},
  \bauthor{\bsnm{Karatzas},~\bfnm{I.}\binits{I.}} \AND
  \bauthor{\bsnm{Yao},~\bfnm{S.}\binits{S.}}
(\byear{2011}).
\bhowpublished{Optimal stopping for dynamic convex risk measures.
  \textit{Illinois J. Math.} (special issue on honor of Don. Burkholder). To
  appear}.
\end{bmisc}
\endbibitem

\bibitem[\protect\citeauthoryear{Bertsekas and Shreve}{1978}]{BerShr78}
\begin{bbook}[mr]
\bauthor{\bsnm{Bertsekas},~\bfnm{Dimitri~P.}\binits{D.~P.}} \AND
  \bauthor{\bsnm{Shreve},~\bfnm{Steven~E.}\binits{S.~E.}}
(\byear{1978}).
\btitle{Stochastic Optimal Control: The Discrete Time Case}.
\bseries{Mathematics in Science and Engineering}
\bvolume{139}.
\bpublisher{Academic Press},
  \baddress{New York}.
\bid{mr={0511544}}
\end{bbook}
\endbibitem

\bibitem[\protect\citeauthoryear{Cheridito, Filipovi{\'c} and
  Yor}{2005}]{CheFilYor05}
\begin{barticle}[mr]
\bauthor{\bsnm{Cheridito},~\bfnm{Patrick}\binits{P.}},
  \bauthor{\bsnm{Filipovi{\'c}},~\bfnm{Damir}\binits{D.}} \AND
  \bauthor{\bsnm{Yor},~\bfnm{Marc}\binits{M.}}
(\byear{2005}).
\btitle{Equivalent and absolutely continuous measure changes for jump-diffusion
  processes}.
\bjournal{Ann. Appl. Probab.}
\bvolume{15}
\bpages{1713--1732}.
\bid{doi={10.1214/105051605000000197}, issn={1050-5164}, mr={2152242}}
\end{barticle}
\endbibitem

\bibitem[\protect\citeauthoryear{Cvitani{\'c}, Pham and
  Touzi}{1999}]{CviPhaTou99}
\begin{barticle}[mr]
\bauthor{\bsnm{Cvitani{\'c}},~\bfnm{Jak{\v{s}}a}\binits{J.}},
  \bauthor{\bsnm{Pham},~\bfnm{Huy{\^e}n}\binits{H.}} \AND
  \bauthor{\bsnm{Touzi},~\bfnm{Nizar}\binits{N.}}
(\byear{1999}).
\btitle{Super-replication in stochastic volatility models under portfolio
  constraints}.
\bjournal{J. Appl. Probab.}
\bvolume{36}
\bpages{523--545}.
\bid{issn={0021-9002}, mr={1724796}}
\end{barticle}
\endbibitem

\bibitem[\protect\citeauthoryear{Delbaen and Schachermayer}{1995a}]{DelSch95N2}
\begin{barticle}[mr]
\bauthor{\bsnm{Delbaen},~\bfnm{Freddy}\binits{F.}} \AND
  \bauthor{\bsnm{Schachermayer},~\bfnm{Walter}\binits{W.}}
(\byear{1995}a).
\btitle{The existence of absolutely continuous local martingale measures}.
\bjournal{Ann. Appl. Probab.}
\bvolume{5}
\bpages{926--945}.
\bid{issn={1050-5164}, mr={1384360}}
\end{barticle}
\endbibitem

\bibitem[\protect\citeauthoryear{Delbaen and Schachermayer}{1995b}]{DelSch95N1}
\begin{barticle}[mr]
\bauthor{\bsnm{Delbaen},~\bfnm{Freddy}\binits{F.}} \AND
  \bauthor{\bsnm{Schachermayer},~\bfnm{Walter}\binits{W.}}
(\byear{1995}b).
\btitle{The no-arbitrage property under a~change of num\'eraire}.
\bjournal{Stochastics Stochastics Rep.}
\bvolume{53}
\bpages{213--226}.
\bid{issn={1045-1129}, mr={1381678}}
\end{barticle}
\endbibitem

\bibitem[\protect\citeauthoryear{Denis and Martini}{2006}]{DenMar06}
\begin{barticle}[mr]
\bauthor{\bsnm{Denis},~\bfnm{Laurent}\binits{L.}} \AND
  \bauthor{\bsnm{Martini},~\bfnm{Claude}\binits{C.}}
(\byear{2006}).
\btitle{A theoretical framework for the pricing of contingent claims in the
  presence of model uncertainty}.
\bjournal{Ann. Appl. Probab.}
\bvolume{16}
\bpages{827--852}.
\bid{doi={10.1214/105051606000000169}, issn={1050-5164}, mr={2244434}}
\end{barticle}
\endbibitem

\bibitem[\protect\citeauthoryear{Ekstr{\"o}m and Tysk}{2004}]{EksTys04}
\begin{bmisc}[mr]
\bauthor{\bsnm{Ekstr{\"o}m},~\bfnm{Erik}\binits{E.}} \AND
  \bauthor{\bsnm{Tysk},~\bfnm{Johan}\binits{J.}}
(\byear{2004}).
\bhowpublished{Comparison of two methods for super-replication. Preprint, Uppsala Univ.}
\end{bmisc}
\endbibitem

\bibitem[\protect\citeauthoryear{Ekstr{\"o}m and Tysk}{2009}]{EksTys09}
\begin{barticle}[mr]
\bauthor{\bsnm{Ekstr{\"o}m},~\bfnm{Erik}\binits{E.}} \AND
  \bauthor{\bsnm{Tysk},~\bfnm{Johan}\binits{J.}}
(\byear{2009}).
\btitle{Bubbles, convexity and the {B}lack--{S}choles equation}.
\bjournal{Ann. Appl. Probab.}
\bvolume{19}
\bpages{1369--1384}.
\bid{doi={10.1214/08-AAP579}, issn={1050-5164}, mr={2538074}}
\end{barticle}
\endbibitem

\bibitem[\protect\citeauthoryear{El~Karoui, H\v{u}\`{u}~Nguyen and
  Jeanblanc-Picqu{\'e}}{1987}]{ElKHuuJea87}
\begin{barticle}[mr]
\bauthor{\bsnm{El~Karoui},~\bfnm{Nicole}\binits{N.}},
\bauthor{\bsnm{H\v{u}\`{u}~Nguyen},~\bfnm{Du'}\binits{D.}} \AND
  \bauthor{\bsnm{Jeanblanc-Picqu{\'e}},~\bfnm{Monique}\binits{M.}}
(\byear{1987}).
\btitle{Compactification methods in the control of degenerate diffusions:
  Existence of an optimal control}.
\bjournal{Stochastics}
\bvolume{20}
\bpages{169--219}.
\bid{issn={0090-9491}, mr={0878312}}
\end{barticle}
\endbibitem

\bibitem[\protect\citeauthoryear{El~Karoui, Jeanblanc-Picqu{\'e} and
  Shreve}{1998}]{ElKJeaShr98}
\begin{barticle}[mr]
\bauthor{\bsnm{El~Karoui},~\bfnm{Nicole}\binits{N.}},
  \bauthor{\bsnm{Jeanblanc-Picqu{\'e}},~\bfnm{Monique}\binits{M.}} \AND
  \bauthor{\bsnm{Shreve},~\bfnm{Steven~E.}\binits{S.~E.}}
(\byear{1998}).
\btitle{Robustness of the {B}lack and {S}choles formula}.
\bjournal{Math. Finance}
\bvolume{8}
\bpages{93--126}.
\bid{doi={10.1111/1467-9965.00047}, issn={0960-1627}, mr={1609962}}
\end{barticle}
\endbibitem

\bibitem[\protect\citeauthoryear{Ethier and Kurtz}{1986}]{EthKur86}
\begin{bbook}[mr]
\bauthor{\bsnm{Ethier},~\bfnm{Stewart~N.}\binits{S.~N.}} \AND
  \bauthor{\bsnm{Kurtz},~\bfnm{Thomas~G.}\binits{T.~G.}}
(\byear{1986}).
\btitle{Markov Processes: Characterization and Convergence}.
\bpublisher{Wiley}, \baddress{New York}.
\bid{doi={10.1002/9780470316658}, mr={0838085}}
\end{bbook}
\endbibitem

\bibitem[\protect\citeauthoryear{Fernholz and Karatzas}{2005}]{FERKAR05}
\begin{barticle}[auto:STB|2011-03-03|12:04:44]
\bauthor{\bsnm{Fernholz},~\bfnm{E.~R.}\binits{E.~R.}} \AND
  \bauthor{\bsnm{Karatzas},~\bfnm{I.}\binits{I.}}
(\byear{2005}).
\btitle{Relative arbitrage in volatility-stabilized markets}.
\bjournal{Annals of Finance}
\bvolume{1}
\bpages{149--177}.
\end{barticle}
\endbibitem

\bibitem[\protect\citeauthoryear{Fernholz and Karatzas}{2009}]{FERKAR}
\begin{bincollection}[auto:STB|2011-03-03|12:04:44]
\bauthor{\bsnm{Fernholz},~\bfnm{E.~R.}\binits{E.~R.}} \AND
  \bauthor{\bsnm{Karatzas},~\bfnm{I.}\binits{I.}}
(\byear{2009}).
\btitle{Stochastic portfolio theory: A survey}.
In \bbooktitle{Handbook of Numerical Analysis}
(\beditor{A. Bensoussan and Q. Zhang}, eds.).
\bpages{88--168}.
\bpublisher{Elsevier}, \baddress{Amsterdam}.
\end{bincollection}
\endbibitem

\bibitem[\protect\citeauthoryear{Fernholz and Karatzas}{2010a}]{FerKar10}
\begin{barticle}[mr]
\bauthor{\bsnm{Fernholz},~\bfnm{Daniel}\binits{D.}} \AND
  \bauthor{\bsnm{Karatzas},~\bfnm{Ioannis}\binits{I.}}
(\byear{2010a}).
\btitle{On optimal arbitrage}.
\bjournal{Ann. Appl. Probab.}
\bvolume{20}
\bpages{1179--1204}.
\bid{doi={10.1214/09-AAP642}, issn={1050-5164}, mr={2676936}}
\end{barticle}
\endbibitem

\bibitem[\protect\citeauthoryear{Fernholz and Karatzas}{2010b}]{FerKar10b}
\begin{bincollection}[auto:STB|2011-03-03|12:04:44]
\bauthor{\bsnm{Fernholz},~\bfnm{Daniel}\binits{D.}} \AND
  \bauthor{\bsnm{Karatzas},~\bfnm{Ioannis}\binits{I.}}
(\byear{2010b}).
\btitle{Probabilistic aspects of arbitrage}.
In
\bbooktitle{Contemporary Mathematical Finance: Essays in Honor of Eckhard Platen}
(\beditor{C. Chiarella and A. Novikov}, eds.)
\bpages{1--17}.
\bpublisher{Springer}, \baddress{New York}.
\end{bincollection}
\endbibitem

\bibitem[\protect\citeauthoryear{Fernholz, Karatzas and
  Kardaras}{2005}]{FERKARKAR05}
\begin{barticle}[auto:STB|2011-03-03|12:04:44]
\bauthor{\bsnm{Fernholz},~\bfnm{E.~R.}\binits{E.~R.}},
  \bauthor{\bsnm{Karatzas},~\bfnm{I.}\binits{I.}} \AND
  \bauthor{\bsnm{Kardaras},~\bfnm{C.}\binits{C.}}
(\byear{2005}).
\btitle{Diversity and arbitrage in equity markets}.
\bjournal{Finance Stoch.}
\bvolume{31}
\bpages{37--53}.
\end{barticle}
\endbibitem

\bibitem[\protect\citeauthoryear{Fleming and Rishel}{1975}]{FleRis75}
\begin{bbook}[mr]
\bauthor{\bsnm{Fleming},~\bfnm{Wendell~H.}\binits{W.~H.}} \AND
  \bauthor{\bsnm{Rishel},~\bfnm{Raymond~W.}\binits{R.~W.}}
(\byear{1975}).
\btitle{Deterministic and Stochastic Optimal Control}.
\bseries{Applications of Mathematics}
\bvolume{1}.
\bpublisher{Springer}, \baddress{Berlin}.
\bid{mr={0454768}}
\end{bbook}
\endbibitem

\bibitem[\protect\citeauthoryear{Fleming and Soner}{1993}]{FleSon93}
\begin{bbook}[mr]
\bauthor{\bsnm{Fleming},~\bfnm{Wendell~H.}\binits{W.~H.}} \AND
  \bauthor{\bsnm{Soner},~\bfnm{H.~Mete}\binits{H.~M.}}
(\byear{1993}).
\btitle{Controlled {M}arkov Processes and Viscosity Solutions}.
\bseries{Applications of Mathematics (New York)}
\bvolume{25}.
\bpublisher{Springer}, \baddress{New York}.
\bid{mr={1199811}}
\end{bbook}
\endbibitem

\bibitem[\protect\citeauthoryear{Fleming and Vermes}{1989}]{FleVer89}
\begin{barticle}[mr]
\bauthor{\bsnm{Fleming},~\bfnm{Wendell~H.}\binits{W.~H.}} \AND
  \bauthor{\bsnm{Vermes},~\bfnm{Domokos}\binits{D.}}
(\byear{1989}).
\btitle{Convex duality approach to the optimal control of diffusions}.
\bjournal{SIAM J. Control Optim.}
\bvolume{27}
\bpages{1136--1155}.
\bid{doi={10.1137/0327060}, issn={0363-0129}, mr={1009341}}
\end{barticle}
\endbibitem

\bibitem[\protect\citeauthoryear{F{\"o}llmer}{1972}]{Fol72}
\begin{barticle}[mr]
\bauthor{\bsnm{F{\"o}llmer},~\bfnm{Hans}\binits{H.}}
(\byear{1972}).
\btitle{The exit measure of a supermartingale}.
\bjournal{Z. Wahrsch. Verw. Gebiete}
\bvolume{21}
\bpages{154--166}.
\bid{mr={0309184}}
\end{barticle}
\endbibitem

\bibitem[\protect\citeauthoryear{F{\"o}llmer}{1973}]{Fol73}
\begin{barticle}[mr]
\bauthor{\bsnm{F{\"o}llmer},~\bfnm{Hans}\binits{H.}}
(\byear{1973}).
\btitle{On the representation of semimartingales}.
\bjournal{Ann. Probab.}
\bvolume{1}
\bpages{580--589}.
\bid{mr={0353446}}
\end{barticle}
\endbibitem

\bibitem[\protect\citeauthoryear{F{\"o}llmer and Gundel}{2006}]{FolGun06}
\begin{barticle}[mr]
\bauthor{\bsnm{F{\"o}llmer},~\bfnm{Hans}\binits{H.}} \AND
  \bauthor{\bsnm{Gundel},~\bfnm{Anne}\binits{A.}}
(\byear{2006}).
\btitle{Robust projections in the class of martingale measures}.
\bjournal{Illinois J. Math.}
\bvolume{50}
\bpages{439--472 (electronic)}.
\bid{issn={0019-2082}, mr={2247836}}
\end{barticle}
\endbibitem

\bibitem[\protect\citeauthoryear{F{\"o}llmer, Schied and
  Weber}{2009}]{FOLSCHWEB}
\begin{bincollection}[auto:STB|2011-03-03|12:04:44]
\bauthor{\bsnm{F{\"o}llmer},~\bfnm{H.}\binits{H.}},
  \bauthor{\bsnm{Schied},~\bfnm{A.}\binits{A.}} \AND
  \bauthor{\bsnm{Weber},~\bfnm{S.}\binits{S.}}
(\byear{2009}).
\btitle{Robust preferences and robust Portfolio choice}.
In \bbooktitle{Handbook of Numerical Analysis}
(\beditor{A. Bensoussan and Q. Zhang}, eds.).
\bpages{29--87}.
\bpublisher{Elsevier}, \baddress{Amsterdam}.
\end{bincollection}
\endbibitem

\bibitem[\protect\citeauthoryear{Frey}{2000}]{Fre00}
\begin{barticle}[mr]
\bauthor{\bsnm{Frey},~\bfnm{R{\"u}diger}\binits{R.}}
(\byear{2000}).
\btitle{Superreplication in stochastic volatility models and optimal stopping}.
\bjournal{Finance Stoch.}
\bvolume{4}
\bpages{161--187}.
\bid{doi={10.1007/s007800050010}, issn={0949-2984}, mr={1780325}}
\end{barticle}
\endbibitem

\bibitem[\protect\citeauthoryear{Gilboa and Schmeidler}{1989}]{GilSch89}
\begin{barticle}[mr]
\bauthor{\bsnm{Gilboa},~\bfnm{Itzhak}\binits{I.}} \AND
  \bauthor{\bsnm{Schmeidler},~\bfnm{David}\binits{D.}}
(\byear{1989}).
\btitle{Maxmin expected utility with nonunique prior}.
\bjournal{J. Math. Econom.}
\bvolume{18}
\bpages{141--153}.
\bid{doi={10.1016/0304-4068(89)90018-9}, issn={0304-4068}, mr={1000102}}
\end{barticle}
\endbibitem

\bibitem[\protect\citeauthoryear{Goia}{2009}]{Goi09}
\begin{bbook}[mr]
\bauthor{\bsnm{Goia},~\bfnm{Irina}\binits{I.}}
(\byear{2009}).
\btitle{Bessel and Volatility-stabilized Processes}.
\bpublisher{ProQuest LLC}, \baddress{Ann Arbor, MI}.
\bnote{PhD Thesis, Columbia Univ}.
\bid{mr={2713615}}
\end{bbook}
\endbibitem

\bibitem[\protect\citeauthoryear{Gozzi and Vargiolu}{2002}]{GozVar02}
\begin{barticle}[mr]
\bauthor{\bsnm{Gozzi},~\bfnm{Fausto}\binits{F.}} \AND
  \bauthor{\bsnm{Vargiolu},~\bfnm{Tiziano}\binits{T.}}
(\byear{2002}).
\btitle{Superreplication of {E}uropean multiasset derivatives with bounded
  stochastic volatility}.
\bjournal{Math. Methods Oper. Res.}
\bvolume{55}
\bpages{69--91}.
\bid{doi={10.1007/s001860200172}, issn={1432-2994}, mr={1892718}}
\end{barticle}
\endbibitem

\bibitem[\protect\citeauthoryear{Gundel}{2005}]{GUN05}
\begin{barticle}[auto:STB|2011-03-03|12:04:44]
\bauthor{\bsnm{Gundel},~\bfnm{A.}\binits{A.}}
(\byear{2005}).
\btitle{Super-replication of European multi-asset derivatives with bounded
  stochastic volatility}.
\bjournal{Finance Stoch.}
\bvolume{9}
\bpages{851--176}.
\end{barticle}
\endbibitem

\bibitem[\protect\citeauthoryear{Haussmann and Lepeltier}{1990}]{HauLep90}
\begin{barticle}[mr]
\bauthor{\bsnm{Haussmann},~\bfnm{U.~G.}\binits{U.~G.}} \AND
  \bauthor{\bsnm{Lepeltier},~\bfnm{J.~P.}\binits{J.~P.}}
(\byear{1990}).
\btitle{On the existence of optimal controls}.
\bjournal{SIAM J. Control Optim.}
\bvolume{28}
\bpages{851--902}.
\bid{doi={10.1137/0328049}, issn={0363-0129}, mr={1051628}}
\end{barticle}
\endbibitem

\bibitem[\protect\citeauthoryear{Heath et~al.}{1987}]{Heaetal87}
\begin{barticle}[mr]
\bauthor{\bsnm{Heath},~\bfnm{D.}\binits{D.}},
  \bauthor{\bsnm{Orey},~\bfnm{S.}\binits{S.}},
  \bauthor{\bsnm{Pestien},~\bfnm{V.}\binits{V.}} \AND
  \bauthor{\bsnm{Sudderth},~\bfnm{W.}\binits{W.}}
(\byear{1987}).
\btitle{Minimizing or maximizing the expected time to reach zero}.
\bjournal{SIAM J. Control Optim.}
\bvolume{25}
\bpages{195--205}.
\bid{doi={10.1137/0325012}, issn={0363-0129}, mr={0872458}}
\end{barticle}
\endbibitem

\bibitem[\protect\citeauthoryear{Janson and Tysk}{2006}]{JanTys06}
\begin{barticle}[mr]
\bauthor{\bsnm{Janson},~\bfnm{Svante}\binits{S.}} \AND
  \bauthor{\bsnm{Tysk},~\bfnm{Johan}\binits{J.}}
(\byear{2006}).
\btitle{Feynman--{K}ac formulas for {B}lack--{S}choles-type operators}.
\bjournal{Bull. London Math. Soc.}
\bvolume{38}
\bpages{269--282}.
\bid{doi={10.1112/S0024609306018194}, issn={0024-6093}, mr={2214479}}
\end{barticle}
\endbibitem

\bibitem[\protect\citeauthoryear{Karatzas and Kardaras}{2007}]{KarKar07}
\begin{barticle}[mr]
\bauthor{\bsnm{Karatzas},~\bfnm{Ioannis}\binits{I.}} \AND
  \bauthor{\bsnm{Kardaras},~\bfnm{Constantinos}\binits{C.}}
(\byear{2007}).
\btitle{The num\'eraire portfolio in semimartingale financial models}.
\bjournal{Finance Stoch.}
\bvolume{11}
\bpages{447--493}.
\bid{doi={10.1007/s00780-007-0047-3}, issn={0949-2984}, mr={2335830}}
\end{barticle}
\endbibitem

\bibitem[\protect\citeauthoryear{Karatzas and Shreve}{1991}]{KarShr91}
\begin{bbook}[mr]
\bauthor{\bsnm{Karatzas},~\bfnm{Ioannis}\binits{I.}} \AND
  \bauthor{\bsnm{Shreve},~\bfnm{Steven~E.}\binits{S.~E.}}
(\byear{1991}).
\btitle{Brownian Motion and Stochastic Calculus},
\bedition{2nd} ed.
\bseries{Graduate Texts in Mathematics}
\bvolume{113}.
\bpublisher{Springer}, \baddress{New York}.
\bid{mr={1121940}}
\end{bbook}
\endbibitem

\bibitem[\protect\citeauthoryear{Karatzas and Zamfirescu}{2005}]{KarZam05}
\begin{barticle}[mr]
\bauthor{\bsnm{Karatzas},~\bfnm{Ioannis}\binits{I.}} \AND
  \bauthor{\bsnm{Zamfirescu},~\bfnm{Ingrid-Mona}\binits{I.-M.}}
(\byear{2005}).
\btitle{Game approach to the optimal stopping problem}.
\bjournal{Stochastics}
\bvolume{77}
\bpages{401--435}.
\bid{doi={10.1080/17442500500219885}, issn={1744-2508}, mr={2178425}}
\end{barticle}
\endbibitem

\bibitem[\protect\citeauthoryear{Kardaras and Robertson}{2011}]{KARROB}
\begin{bmisc}[auto:STB|2011-03-03|12:04:44]
\bauthor{\bsnm{Kardaras},~\bfnm{C.}\binits{C.}} \AND
  \bauthor{\bsnm{Robertson},~\bfnm{S.}\binits{S.}}
(\byear{2011}).
\bhowpublished{Robust maximization of asymptotic growth. Preprint,
  Boston Univ.}
\end{bmisc}
\endbibitem

\bibitem[\protect\citeauthoryear{Krylov}{1973}]{Kry73}
\begin{barticle}[mr]
\bauthor{\bsnm{Krylov},~\bfnm{N.~V.}\binits{N.~V.}}
(\byear{1973}).
\btitle{The selection of a {M}arkov process from a {M}arkov system of
  processes, and the construction of quasidiffusion processes}.
\bjournal{Izv. Akad. Nauk SSSR Ser. Mat.}
\bvolume{37}
\bpages{691--708}.
\bid{issn={0373-2436}, mr={0339338}}
\end{barticle}
\endbibitem

\bibitem[\protect\citeauthoryear{Krylov}{1980}]{Kry80}
\begin{bbook}[mr]
\bauthor{\bsnm{Krylov},~\bfnm{N.~V.}\binits{N.~V.}}
(\byear{1980}).
\btitle{Controlled Diffusion Processes}.
\bseries{Applications of Mathematics}
\bvolume{14}.
\bpublisher{Springer}, \baddress{New York}.
\bid{mr={0601776}}
\end{bbook}
\endbibitem

\bibitem[\protect\citeauthoryear{Krylov}{1987}]{Kry87}
\begin{bbook}[mr]
\bauthor{\bsnm{Krylov},~\bfnm{N.~V.}\binits{N.~V.}}
(\byear{1987}).
\btitle{Nonlinear Elliptic and Parabolic Equations of the Second Order}.
\bseries{Mathematics and Its Applications (Soviet Series)}
\bvolume{7}.
\bpublisher{Reidel}, \baddress{Dordrecht}.
\bid{mr={0901759}}
\end{bbook}
\endbibitem


\bibitem[\protect\citeauthoryear{Krylov}{1989}]{Kry89N1}
\begin{barticle}[mr]
\bauthor{\bsnm{Krylov},~\bfnm{N.~V.}\binits{N.~V.}}
(\byear{1989}).
\btitle{A supermartingale characterization of a set of stochastic integrals}.
\bjournal{Ukrain. Mat. Zh.}
\bvolume{41}
\bpages{757--762, 861}.
\bid{doi={10.1007/BF01060562}, issn={0041-6053}, mr={1002711}}
\end{barticle}
\endbibitem

\bibitem[\protect\citeauthoryear{Krylov}{1990}]{Kry89N2}
\begin{barticle}[mr]
\bauthor{\bsnm{Krylov},~\bfnm{N.~V.}\binits{N.~V.}}
(\byear{1990}).
\btitle{Smoothness of the value function for a controlled diffusion process in a domain}.
\bjournal{Math. USSR Izvestiya}
\bvolume{34}
\bpages{65--95}.
\end{barticle}
\endbibitem

\bibitem[\protect\citeauthoryear{Krylov}{2002}]{Kry02}
\begin{barticle}[mr]
\bauthor{\bsnm{Krylov},~\bfnm{N.~V.}\binits{N.~V.}}
(\byear{2002}).
\btitle{A supermartingale characterization of sets of stochastic integrals and
  applications}.
\bjournal{Probab. Theory Related Fields}
\bvolume{123}
\bpages{521--552}.
\bid{doi={10.1007/s004400100190}, issn={0178-8051}, mr={1921012}}
\end{barticle}
\endbibitem

\bibitem[\protect\citeauthoryear{Kunita}{1990}]{Kun90}
\begin{bbook}[mr]
\bauthor{\bsnm{Kunita},~\bfnm{Hiroshi}\binits{H.}}
(\byear{1990}).
\btitle{Stochastic Flows and Stochastic Differential Equations}.
\bseries{Cambridge Studies in Advanced Mathematics}
\bvolume{24}.
\bpublisher{Cambridge Univ. Press}, \baddress{Cambridge}.
\bid{mr={1070361}}
\end{bbook}
\endbibitem

\bibitem[\protect\citeauthoryear{Levental and Skorohod}{1995}]{LevSko95}
\begin{barticle}[mr]
\bauthor{\bsnm{Levental},~\bfnm{Shlomo}\binits{S.}} \AND
  \bauthor{\bsnm{Skorohod},~\bfnm{Anatolii~V.}\binits{A.~V.}}
(\byear{1995}).
\btitle{A necessary and sufficient condition for absence of arbitrage with tame
  portfolios}.
\bjournal{Ann. Appl. Probab.}
\bvolume{5}
\bpages{906--925}.
\bid{issn={1050-5164}, mr={1384359}}
\end{barticle}
\endbibitem

\bibitem[\protect\citeauthoryear{Lieberman}{1996}]{Lie96}
\begin{bbook}[mr]
\bauthor{\bsnm{Lieberman},~\bfnm{Gary~M.}\binits{G.~M.}}
(\byear{1996}).
\btitle{Second Order Parabolic Differential Equations}.
\bpublisher{World Scientific}, \baddress{River Edge, NJ}.
\bid{mr={1465184}}
\end{bbook}
\endbibitem

\bibitem[\protect\citeauthoryear{Lions}{1983a}]{Lio83N1}
\begin{barticle}[mr]
\bauthor{\bsnm{Lions},~\bfnm{P.~L.}\binits{P.~L.}}
(\byear{1983}a).
\btitle{Optimal control of diffusion processes and
  {H}amilton--{J}acobi--{B}ellman equations. {I}. {T}he dynamic programming
  principle and applications}.
\bjournal{Comm. Partial Differential Equations}
\bvolume{8}
\bpages{1101--1174}.
\bid{doi={10.1080/03605308308820297}, issn={0360-5302}, mr={0709164}}
\end{barticle}
\endbibitem

\bibitem[\protect\citeauthoryear{Lions}{1983b}]{Lio83N2}
\begin{barticle}[mr]
\bauthor{\bsnm{Lions},~\bfnm{P.~L.}\binits{P.~L.}}
(\byear{1983}b).
\btitle{Optimal control of diffusion processes and
  {H}amilton--{J}acobi--{B}ellman equations. {II}. {V}iscosity solutions and
  uniqueness}.
\bjournal{Comm. Partial Differential Equations}
\bvolume{8}
\bpages{1229--1276}.
\bid{doi={10.1080/03605308308820301}, issn={0360-5302}, mr={0709162}}
\end{barticle}
\endbibitem

\bibitem[\protect\citeauthoryear{Lions}{1983c}]{Lio83N3}
\begin{bincollection}[mr]
\bauthor{\bsnm{Lions},~\bfnm{P.~L.}\binits{P.~L.}}
(\byear{1983}c).
\btitle{Optimal control of diffusion processes and
  {H}amilton--{J}acobi--{B}ellman equations. {III}. {R}egularity of the optimal
  cost function}.
In \bbooktitle{Nonlinear Partial Differential Equations and Their Applications.
  {C}oll\`ege de {F}rance Seminar, {V}ol. {V} ({P}aris, 1981/1982)}.
\bseries{Res. Notes in Math.}
\bvolume{93}
\bpages{95--205}.
\bpublisher{Pitman}, \baddress{Boston, MA}.
\bid{mr={0725360}}
\end{bincollection}
\endbibitem

\bibitem[\protect\citeauthoryear{Lions}{1984}]{Lio84}
\begin{bincollection}[mr]
\bauthor{\bsnm{Lions},~\bfnm{Pierre-Louis}\binits{P.-L.}}
(\byear{1984}).
\btitle{Some recent results in the optimal control of diffusion processes}.
In \bbooktitle{Stochastic Analysis ({K}atata/{K}yoto, 1982)}.
\bseries{North-Holland Math. Library}
\bvolume{32}
\bpages{333--367}.
\bpublisher{North-Holland}, \baddress{Amsterdam}.
\bid{mr={0780764}}
\bptnote{check year}
\end{bincollection}
\endbibitem

\bibitem[\protect\citeauthoryear{Lyons}{1995}]{LYO95}
\begin{barticle}[auto:STB|2011-03-03|12:04:44]
\bauthor{\bsnm{Lyons},~\bfnm{T.~J.}\binits{T.~J.}}
(\byear{1995}).
\btitle{Uncertain volatility and the risk-free synthesis of securities}.
\bjournal{Appl. Math. Finance}
\bvolume{2}
\bpages{117--133}.
\end{barticle}
\endbibitem

\bibitem[\protect\citeauthoryear{Meyer}{1972}]{Mey72}
\begin{bincollection}[mr]
\bauthor{\bsnm{Meyer},~\bfnm{P.~A.}\binits{P.~A.}}
(\byear{1972}).
\btitle{La mesure de {H}. {F}\"ollmer en th\'eorie des surmartingales}.
In \bbooktitle{S\'eminaire de {P}robabilit\'es, {VI} ({U}niv. {S}trasbourg,
  Ann\'ee Universitaire 1970--1971; {J}ourn\'ees {P}robabilistes de
  {S}trasbourg, 1971)}.
\bseries{Lecture Notes in Math.}
\bvolume{258}
\bpages{118--129}.
\bpublisher{Springer}, \baddress{Berlin}.
\bid{mr={0368131}}
\end{bincollection}
\endbibitem

\bibitem[\protect\citeauthoryear{Meyer}{2006}]{Mey06}
\begin{barticle}[mr]
\bauthor{\bsnm{Meyer},~\bfnm{Gunter~H.}\binits{G.~H.}}
(\byear{2006}).
\btitle{The {B}lack {S}choles {B}arenblatt equation for options with uncertain
  volatility and its application to static hedging}.
\bjournal{Int. J. Theor. Appl. Finance}
\bvolume{9}
\bpages{673--703}.
\bid{doi={10.1142/S0219024906003755}, issn={0219-0249}, mr={2254127}}
\end{barticle}
\endbibitem

\bibitem[\protect\citeauthoryear{Nutz}{2010}]{NUT}
\begin{bmisc}[auto:STB|2011-03-03|12:04:44]
\bauthor{\bsnm{Nutz},~\bfnm{M.}\binits{M.}}
(\byear{2010}).
\bhowpublished{Random G-expectations. Preprint, ETH Z\"urich.}
\end{bmisc}
\endbibitem

\bibitem[\protect\citeauthoryear{Orey, Pestien and
  Sudderth}{1987}]{OrePesSud87}
\begin{barticle}[mr]
\bauthor{\bsnm{Orey},~\bfnm{Steven}\binits{S.}},
  \bauthor{\bsnm{Pestien},~\bfnm{Victor}\binits{V.}} \AND
  \bauthor{\bsnm{Sudderth},~\bfnm{William}\binits{W.}}
(\byear{1987}).
\btitle{Reaching zero rapidly}.
\bjournal{SIAM J. Control Optim.}
\bvolume{25}
\bpages{1253--1265}.
\bid{doi={10.1137/0325069}, issn={0363-0129}, mr={0905044}}
\end{barticle}
\endbibitem

\bibitem[\protect\citeauthoryear{Pal}{2011}]{PAL}
\begin{bmisc}[auto:STB|2011-03-03|12:04:44]
\bauthor{\bsnm{Pal},~\bfnm{S.}\binits{S.}}
(\byear{2011}).
\bhowpublished{Analysis of the market weights under the volatility-stabilized
  market models. \textit{Ann. Appl. Probab.} To appear}.
\end{bmisc}
\endbibitem

\bibitem[\protect\citeauthoryear{Pal and Protter}{2010}]{PALPRO10}
\begin{barticle}[auto:STB|2011-03-03|12:04:44]
\bauthor{\bsnm{Pal},~\bfnm{S.}\binits{S.}} \AND
  \bauthor{\bsnm{Protter},~\bfnm{Ph.}\binits{P.}}
(\byear{2010}).
\btitle{Analysis of continuous strict local martingales via $h$-transforms}.
\bjournal{Stochastic Process. Appl.}
\bvolume{120}
\bpages{1424--1443}.
\end{barticle}
\endbibitem

\bibitem[\protect\citeauthoryear{Parthasarathy}{1967}]{Par67}
\begin{bbook}[mr]
\bauthor{\bsnm{Parthasarathy},~\bfnm{K.~R.}\binits{K.~R.}}
(\byear{1967}).
\btitle{Probability Measures on Metric Spaces}.
\bseries{Probability and Mathematical Statistics}
\bvolume{3}.
\bpublisher{Academic Press}, \baddress{New York}.
\bid{mr={0226684}}
\end{bbook}
\endbibitem

\bibitem[\protect\citeauthoryear{Peng}{2010}]{PEN}
\begin{bmisc}[auto:STB|2011-03-03|12:04:44]
\bauthor{\bsnm{Peng},~\bfnm{S.}\binits{S.}}
(\byear{2010}).
\bhowpublished{Nonlinear expectations and stochastic calculus under
  uncertainty. Preprint, Shandong Univ.}
\end{bmisc}
\endbibitem

\bibitem[\protect\citeauthoryear{Pestien and Sudderth}{1985}]{PesSud85}
\begin{barticle}[mr]
\bauthor{\bsnm{Pestien},~\bfnm{Victor~C.}\binits{V.~C.}} \AND
  \bauthor{\bsnm{Sudderth},~\bfnm{William~D.}\binits{W.~D.}}
(\byear{1985}).
\btitle{Continuous-time red and black: How to control a diffusion to a goal}.
\bjournal{Math. Oper. Res.}
\bvolume{10}
\bpages{599--611}.
\bid{doi={10.1287/moor.10.4.599}, issn={0364-765X}, mr={0812818}}
\end{barticle}
\endbibitem

\bibitem[\protect\citeauthoryear{Protter}{2004}]{Pro04}
\begin{bbook}[mr]
\bauthor{\bsnm{Protter},~\bfnm{Philip~E.}\binits{P.~E.}}
(\byear{2004}).
\btitle{Stochastic Integration and Differential Equations},
\bedition{2nd} ed.
\bseries{Applications of Mathematics (New York)}
\bvolume{21}.
\bpublisher{Springer}, \baddress{Berlin}.
\bid{mr={2020294}}
\end{bbook}
\endbibitem

\bibitem[\protect\citeauthoryear{Riedel}{2009}]{Rie09}
\begin{barticle}[mr]
\bauthor{\bsnm{Riedel},~\bfnm{Frank}\binits{F.}}
(\byear{2009}).
\btitle{Optimal stopping with multiple priors}.
\bjournal{Econometrica}
\bvolume{77}
\bpages{857--908}.
\bid{doi={10.3982/ECTA7594}, issn={0012-9682}, mr={2531363}}
\end{barticle}
\endbibitem

\bibitem[\protect\citeauthoryear{Romagnoli and Vargiolu}{2000}]{RomVar00}
\begin{barticle}[mr]
\bauthor{\bsnm{Romagnoli},~\bfnm{Silvia}\binits{S.}} \AND
  \bauthor{\bsnm{Vargiolu},~\bfnm{Tiziano}\binits{T.}}
(\byear{2000}).
\btitle{Robustness of the {B}lack--{S}choles approach in the case of options on
  several assets}.
\bjournal{Finance Stoch.}
\bvolume{4}
\bpages{325--341}.
\bid{doi={10.1007/s007800050076}, issn={0949-2984}, mr={1779582}}
\end{barticle}
\endbibitem

\bibitem[\protect\citeauthoryear{Ruf}{2011}]{RUF}
\begin{bmisc}[auto:STB|2011-03-03|12:04:44]
\bauthor{\bsnm{Ruf},~\bfnm{J.}\binits{J.}}
(\byear{2011}).
\bhowpublished{Hedging under arbitrage. \textit{Math. Finance}.
To appear}.
\end{bmisc}
\endbibitem

\bibitem[\protect\citeauthoryear{Schied}{2007}]{Sch07}
\begin{barticle}[mr]
\bauthor{\bsnm{Schied},~\bfnm{Alexander}\binits{A.}}
(\byear{2007}).
\btitle{Optimal investments for risk- and ambiguity-averse preferences:
A~duality approach}.
\bjournal{Finance Stoch.}
\bvolume{11}
\bpages{107--129}.
\bid{doi={10.1007/s00780-006-0024-2}, issn={0949-2984}, mr={2284014}}
\end{barticle}
\endbibitem

\bibitem[\protect\citeauthoryear{Schied and Wu}{2005}]{SchWu05}
\begin{barticle}[mr]
\bauthor{\bsnm{Schied},~\bfnm{Alexander}\binits{A.}} \AND
  \bauthor{\bsnm{Wu},~\bfnm{Ching-Tang}\binits{C.-T.}}
(\byear{2005}).
\btitle{Duality theory for optimal investments under model uncertainty}.
\bjournal{Statist. Decisions}
\bvolume{23}
\bpages{199--217}.
\bid{doi={10.1524/stnd.2005.23.3.199}, issn={0721-2631}, mr={2236457}}
\end{barticle}
\endbibitem

\bibitem[\protect\citeauthoryear{Soner, Touzi and Zhang}{2010a}]{SONTOUZHAN1}
\begin{bmisc}[auto:STB|2011-03-03|12:04:44]
\bauthor{\bsnm{Soner},~\bfnm{H.~M.}\binits{H.~M.}},
  \bauthor{\bsnm{Touzi},~\bfnm{N.}\binits{N.}} \AND
  \bauthor{\bsnm{Zhang},~\bfnm{J.}\binits{J.}}
(\byear{2010}a).
\bhowpublished{Quasi-sure stochastic analysis through aggregation. Preprint, ETH
  Z\"urich.}
\end{bmisc}
\endbibitem

\bibitem[\protect\citeauthoryear{Soner, Touzi and Zhang}{2010b}]{SONTOUZHAN2}
\begin{bmisc}[auto:STB|2011-03-03|12:04:44]
\bauthor{\bsnm{Soner},~\bfnm{H.~M.}\binits{H.~M.}},
  \bauthor{\bsnm{Touzi},~\bfnm{N.}\binits{N.}} \AND
  \bauthor{\bsnm{Zhang},~\bfnm{J.}\binits{J.}}
(\byear{2010}b).
\bhowpublished{Dual formulation of second-order target problems.  Preprint, ETH Z\"urich.
}
\end{bmisc}
\endbibitem

\bibitem[\protect\citeauthoryear{Stroock and Varadhan}{1979}]{StrVar79}
\begin{bbook}[mr]
\bauthor{\bsnm{Stroock},~\bfnm{Daniel~W.}\binits{D.~W.}} \AND
  \bauthor{\bsnm{Varadhan},~\bfnm{S.~R.~Srinivasa}\binits{S.~R.~S.}}
(\byear{1979}).
\btitle{Multidimensional Diffusion Processes}.
\bseries{Grundlehren der Mathematischen Wissenschaften [Fundamental Principles
  of Mathematical Sciences]}
\bvolume{233}.
\bpublisher{Springer}, \baddress{Berlin}.
\bid{mr={0532498}}
\end{bbook}
\endbibitem

\bibitem[\protect\citeauthoryear{Sudderth and Weerasinghe}{1989}]{SudWee89}
\begin{barticle}[mr]
\bauthor{\bsnm{Sudderth},~\bfnm{William~D.}\binits{W.~D.}} \AND
  \bauthor{\bsnm{Weerasinghe},~\bfnm{Ananda}\binits{A.}}
(\byear{1989}).
\btitle{Controlling a process to a goal in finite time}.
\bjournal{Math. Oper. Res.}
\bvolume{14}
\bpages{400--409}.
\bid{doi={10.1287/moor.14.3.400}, issn={0364-765X}, mr={1008420}}
\end{barticle}
\endbibitem

\bibitem[\protect\citeauthoryear{Talay and Zheng}{2002}]{TalZhe02}
\begin{barticle}[mr]
\bauthor{\bsnm{Talay},~\bfnm{Denis}\binits{D.}} \AND
  \bauthor{\bsnm{Zheng},~\bfnm{Ziyu}\binits{Z.}}
(\byear{2002}).
\btitle{Worst case model risk management}.
\bjournal{Finance Stoch.}
\bvolume{6}
\bpages{517--537}.
\bid{doi={10.1007/s007800200074}, issn={0949-2984}, mr={1932383}}
\end{barticle}
\endbibitem

\bibitem[\protect\citeauthoryear{Van~Mellaert and Dorato}{1972}]{VanDOR72}
\begin{barticle}[auto:STB|2011-03-03|12:04:44]
\bauthor{\bsnm{Van~Mellaert},~\bfnm{L.~J.}\binits{L.~J.}} \AND
  \bauthor{\bsnm{Dorato},~\bfnm{P.}\binits{P.}}
(\byear{1972}).
\btitle{Numerical solution of an optimal control problem with a probability
  criterion}.
\bjournal{IEEE Transactions on Automatic Control}
\bvolume{AC-17}
\bpages{543--546}.
\end{barticle}
\endbibitem

\bibitem[\protect\citeauthoryear{Vargiolu}{2001}]{VAR}
\begin{bmisc}[auto:STB|2011-03-03|12:04:44]
\bauthor{\bsnm{Vargiolu},~\bfnm{T.}\binits{T.}}
(\byear{2001}).
\bhowpublished{Existence, uniqueness and smoothness for the
  Black--Scholes--Barenblatt equation. Technical report, Univ. Padova.}
\end{bmisc}
\endbibitem

\bibitem[\protect\citeauthoryear{Vorbrink}{2010}]{VOR}
\begin{bmisc}[auto:STB|2011-03-03|12:04:44]
\bauthor{\bsnm{Vorbrink},~\bfnm{J.}\binits{J.}}
(\byear{2010}).
\bhowpublished{Financial markets with volatility uncertainty. Technical
report,
  Univ. Bielefeld.}
\end{bmisc}
\endbibitem

\bibitem[\protect\citeauthoryear{Wang}{1992a}]{WAN92}
\begin{barticle}[mr]
\bauthor{\bsnm{Wang},~\bfnm{Lihe}\binits{L.}}
(\byear{1992}a).
\btitle{On the regularity theory of fully nonlinear parabolic equations. {I}}.
\bjournal{Comm. Pure Appl. Math.}
\bvolume{45}
\bpages{27--76}.
\bid{doi={10.1002/cpa.3160450103}, issn={0010-3640}, mr={1135923}}
\end{barticle}
\endbibitem

\bibitem[\protect\citeauthoryear{Wang}{1992b}]{Wan92N2}
\begin{barticle}[mr]
\bauthor{\bsnm{Wang},~\bfnm{Lihe}\binits{L.}}
(\byear{1992}b).
\btitle{On the regularity theory of fully nonlinear parabolic equations. {II}}.
\bjournal{Comm. Pure Appl. Math.}
\bvolume{45}
\bpages{141--178}.
\bid{doi={10.1002/cpa.3160450202}, issn={0010-3640}, mr={1139064}}
\end{barticle}
\endbibitem

\bibitem[\protect\citeauthoryear{Wang}{1992c}]{Wan92N3}
\begin{barticle}[mr]
\bauthor{\bsnm{Wang},~\bfnm{Lihe}\binits{L.}}
(\byear{1992}c).
\btitle{On the regularity theory of fully nonlinear parabolic equations.
  {III}}.
\bjournal{Comm. Pure Appl. Math.}
\bvolume{45}
\bpages{255--262}.
\bid{doi={10.1002/cpa.3160450302}, issn={0010-3640}, mr={1151267}}
\end{barticle}
\endbibitem

\end{thebibliography}
\end{document}